\documentclass[11pt,a4paper]{article}
\pdfoutput=1
\usepackage{color,graphicx}
\usepackage[utf8]{inputenc}
\usepackage[T1]{fontenc}
\usepackage{amsmath,amsbsy,amsfonts,amssymb}
\usepackage{theorem,wasysym}
\usepackage{graphicx,psfrag}
\usepackage{dsfont}
\usepackage{hyperref}
\usepackage{algorithm}
\usepackage{algorithmic}
\usepackage{natbib}
\usepackage{subcaption}
\usepackage{lineno}
\usepackage{setspace}
\usepackage{geometry}
\def\be{\beta}

\def\Ga{\Gamma}

\def\eps{\varepsilon}
\def\tet{\theta}
\def\Tet{\Theta}
\def\lam{\lambda}
\def\Lam{\Lambda}

\def\sig{\sigma}
\def\Sig{\Sigma}

\def\f{\times}

\def\E{\mathds{E}}
\def\V{\mathds{V}}
\def\P{\mathds{P}}

\def\R{\mathds{R}}

\def\No{\mathcal{N}}
\def\X{\mathcal{X}}

\def\Cov{\mathrm{Cov}}
\def\A{\mathrm{A}}

\def\C{\mathrm{C}}

\def\new{\mathrm{new}}

\newcommand{\tr}{\mbox{Tr}}
\begin{document}

\title{Combining covariance tapering and lasso driven low rank
  decomposition for the kriging of large spatial datasets}


\author{Thomas Romary \and Nicolas Desassis}




\maketitle

\begin{abstract}
  Large spatial datasets are becoming ubiquitous in environmental
  sciences with the explosion in the amount of data produced by
  sensors that monitor and measure the Earth system. Consequently, the
  geostatistical analysis of these data requires adequate methods.
  Richer datasets lead to more complex modeling but may also prevent
  from using classical techniques. Indeed, the kriging predictor is
  not straightforwarldly available as it requires the inversion of the
  covariance matrix of the data. The challenge of handling such
  datasets is therefore to extract the maximum of information they
  contain while ensuring the numerical tractability of the associated
  inference and prediction algorithms. The different approaches that
  have been developed in the literature to address this problem can be
  classified into two families, both aiming at making the inversion of
  the covariance matrix computationally feasible. The covariance
  tapering approach circumvents the problem by enforcing the sparsity
  of the covariance matrix, making it invertible in a reasonable
  computation time. The second available approach assumes a low rank
  representation of the covariance function. While both approaches
  have their drawbacks, we propose a way to combine them and benefit
  from their advantages. The covariance model is assumed to have the
  form low rank plus sparse. The choice of the basis functions
  sustaining the low rank component is data driven and is achieved
  through a selection procedure, thus alleviating the computational
  burden of the low rank part. This model expresses as a spatial
  random effects model and the estimation of the parameters is
  conducted through a step by step approach treating each scale
  separately. The resulting model can account for second order non
  stationarity and handle large volumes of data.

\end{abstract}

\section{Introduction}
\label{sec:intro}
While a spatial datum was expensive to obtain in the traditional
application fields of Geostatistics (e.g. drilling wells for oil
reserve estimation), with the increasing deployment of remote sensing
platforms and sensors networks, spatial data-base paradigms have moved
from small to massive, often of the order of gigabytes per
day. Therefore, methods for the geostatistical analysis of these kinds
of data have to be developed. Indeed, richer datasets allows for more
complex modeling but may also prevent from using straightforwardly
classical techniques. The challenge of handling such datasets is to
extract the maximum of information they contain while ensuring the
numerical tractability of the associated inference and prediction
algorithms. As an example, satellite image restoration is a
challenging problem for geostatisticians, as it involves large amounts
of data and possible non stationarity over the domain of interest,
both in space and time.

The classical spatial predictor in geostatistics is the kriging, see
e.g. \citep{fasc5,chiles}. It applies when the phenomenon under study can
be modeled by a second order random field. The first order moment
accounts generally for large scale fluctuations of the phenomenon. It
is described in a basis of known functions of the space-time,
e.g. polynomials or exhaustively known covariates. They are generally
called the drift functions. The covariance of the random field
describes the regularity of the medium to small scale variability of
the phenomenon.

The inference is generally conducted as follows: the drift functions
are chosen. Then a parametric covariance model is selected and
fitted. This step can also be straightforwardly automatized through
the use of model selection and fitting algorithms. This step is
conducted either through a method of moments
\citep{desassis2013automatic} or maximum-likelihood under a Gaussian
hypothesis \citep{mardia}. This inference step is made more complex
for several reasons in the context of large datasets. The phenomenon
may exhibit a nonstationary covariance structure that has to be
inferred \citep{Fouedjio201545,fouedjio2016generalized}. The
computation of the likelihood is also made impossible as it involves
the computation of the determinant and the inverse of the covariance
matrix of the data, whose size is $n \f n$, where $n$ is the number of
data. An approximation of the likelihood by compositing can then be
necessary \citep{stein2004approximating,
  varin2011overview,eidsvik2014estimation}.

Once the covariance model is known, the kriging interpolator is built
as the best linear unbiased predictor of the random field at a new
location. It coincides with the conditional expectation under Gaussian
hypothesis. It consists of a weighted average of the data. The
computation of the weights involves the inversion of the covariance
matrix, whose size is $n\f n$. This approach is therefore intractable
when $n$ becomes large, as the computational complexity of the
inversion scales with $n^3$.

Several approaches have been developed in the literature to address
this problem. These approaches can be classified into two families,
both aiming at making the inversion of the covariance matrix
computationally feasible. 

The first one circumvents the problem by enforcing the sparsity of the
covariance matrix \citep{tapering}, making it invertible in a
reasonable computation time ($O(n)$). This implies however that the
covariance model that is fitted to the data is compactly supported
with a short range \citep{compcov}, neglecting possible medium scale
variability. In other words, the fitted covariance is constrained to
this particular structure and may not reflect the whole complexity of
the phenomenon fluctuations. More precisely, while accurate in densely
sampled areas, this method generally fails at providing correct
predictions in scarcely sampled areas, as the correlation between
remote locations decreases quickly towards zero, especially at the
edges of the domain \citep{stein2013statistical}. This approach
simplifies the inference. In particular, it makes the
maximum-likelihood approach tractable as the inverse covariance matrix
is made accessible \citep{taperingest}. 

The second available approach assumes, on the contrary, a low rank
representation of the covariance function. Many low rank approaches
have been proposed over the years, among which we can cite the fixed
rank kriging \citep{fixedrank}, the predictive processes
\citep{pred_proc} and also more recent works
\citep{katzfuss2017multi,nychka2015multiresolution,incholkrig}. It
assumes that the covariance function can be represented on a limited
number, say $p$, of basis functions. Adding a nugget effect, i.e. a
local unstructured noise, to the modeling makes the inversion of the
covariance matrix feasible through the Woodbury matrix inversion
formula with a complexity $O(n \f p^3)$, where $p$ is small,
relatively to $n$. The latter approach tends to overestimate the
regularity of the spatial phenomenon, failing to capture its local
behavior. Again, the whole complexity of the phenomenon fluctuations
may not be correctly reproduced in this framework, as the small scale
fluctuatione are considered unstructured and modelled by a nugget
effect. The generated predictions are generally too smooth, failing at
reproducing the small scale fluctuations, see
\citep{stein2014limitations}. It is worth noting that the low rank
approach provides a non stationary covariance model by
construction. Concerning the inference, chiefly two approaches are
available: either a conventional model is fitted and a low rank
approximation is consequently computed \citep{pred_proc,incholkrig},
or the inference is conducted through an expectation-maximization (EM)
algorithm \citep{katzfuss2011spatio}. Direct method of moments
approach is also available \citep{fixedrank}.

Finally, Markov random field models depend on the observation
locations, and realignment to a much finer grid with missing values is
required for irregular locations \citep{spde,sun2016statistically},
which induces an additional cost for both inference and prediction.

In this work, we propose to use covariance models that combines the
advantages of both approaches. Indeed, using additively a compactly
supported short range covariance and a low rank component in the
modeling of the covariance allows capturing both local and medium
scale fluctuations of the phenomenon and therefore improving the
prediction accuracy. In this settings, the covariance matrix is still
invertible through the Woodbury formula at the same computational cost
as the low rank approach ($O(n \f p^3)$). This is not a new idea.  It
has been already proposed in \citep{stein_kor}, where a few Legendre
polynomials are used on the sphere, in \citep{full_scale_cov} where
the remainder of a predictive process approximation is tapered and
more recently in \citep{ma2017fused}. The originality of our work is
to take advantage of a scale separation between the sparse and low
rank terms. Consequently, the inference of both terms can be conducted
almost separately, making it computationally efficient.

Further more, in low rank approximations, the choice of the basis
functions is generally taken arbitrarily. The set of basis functions
is either built from a fixed family of functions
\citep{fixedrank,wnychka,stein_kor} or from colums of stationary
isotropic covariance matrices
\citep{pred_proc,katzfuss2017multi}. Here, the choice of the basis
functions sustaining the low rank component is data driven and
achieved through a selection procedure, which results in practice in a
low number of functions, hence decreasing the associated computational
burden. Moreover, this provides an important flexibility to the model
allowing it to capture a large variety of possible fluctuations
including varying anisotropy and smoothness, as illustrated in the
examples.

The paper is organized as follows. The first part reviews briefly
existing approaches and describes the proposed one. The second part
details the developed inference procedure. Finally, the performances
of the approach are assessed on synthetic examples and a real dataset.

\section{Modeling}
\label{sec:mod}
We assume that the phenomenon under study can be modeled by a second
order random field of the form:
\begin{equation}
  Z(x) = \mu(x) + Y(x),~x\in\X\subset\R^d,
\end{equation}
where $\mu(x)$ is a deterministic trend and $Y(x)$ is a centered
second order random field with covariance function
$\Cov(Y(x),Y(y)) = \C(x,y)$. More specifically, the trend can be
expressed as a weighted sum of drift functions:
\begin{equation}
  \mu(x) = \sum_{l=1}^L\be_lf_l(x),
\end{equation}
where $(\be_l)_{l=1,\ldots,L}$ and $(f_l)_{l=1,\ldots,L}$ are
deterministic and $L$ is fixed. The drift functions are functions of
the space, such as polynomials, and/or auxiliary variables accounting
for a systematic behaviour of the phenomenon under study. The first
drift function is generally the indicator of $\X$, accounting for a
constant mean. A selection of the relevant functions can be performed
preliminarily using standard variable selection techniques such as
forward, backward, stepwise methods \citep[see e.g.][]{saporta} or
penalization approaches \citep[see e.g.][]{elements} if the candidates
are many.  This modelling is common in Geostatistics where the trend
accounts for the large scale fluctuations while the second term
accounts for smaller scale fluctuations occurring around the drift.

\subsection{Reminder on kriging}
Under these settings, given a set of observations $z =
(z(x_1),\ldots,z(x_n))^t$ from a realization of $Z$ sampled at
locations $(x_1,\ldots,x_n)\in\X^n$, the universal kriging predictor
provides the best linear unbiased predictor in the sense of the mean
squared error at a new location $x_0$. It is a weighted sum of the
observations:
\begin{equation}
  Z^* = \Lam^tz,
\end{equation}
whose weights $\Lam = (\lam_1,\ldots,\lam_n)$ are solution of the
universal kriging system, straightforwardly obtained minimizing the
mean squared prediction error under unbiasedness constraints:
\begin{equation}
\label{eq:kriging}  
\begin{pmatrix} C & F\\ F^t &
    0 \end{pmatrix} \begin{pmatrix} \Lam \\ \mu\end{pmatrix}
    = \begin{pmatrix} C_0 \\ f_0\end{pmatrix}
\end{equation}
where $C$ is the covariance matrix associated to the observations,
i.e. with entries $C_{ij} = \C(x_i,x_j),~(i,j)\in\{1,\ldots,n\}^2$,
$C_0$ is the vector of covariances between $Z(x_0)$ and the
observations, $f_0= (f_1(x_0),\ldots,f_l(x_0))^t$, $$F
= \begin{pmatrix} f_1(x_1) & \ldots & f_L(x_1) \\ \vdots & \ddots &
  \vdots \\ f_1(x_n) & \ldots & f_L(x_n)\end{pmatrix},$$ and $\mu =
(\mu_1,\ldots,\mu_L)^t$ is the vector of Lagrange multipliers
accounting for the unbiasedness constraints.

The solution to this system is obtained by block inversion:
\begin{equation}
  \label{eq:uk}
  Z^*=C_0^tC^{-1}z + (f_0^t-C_0^tC^{-1}F)(F^tC^{-1}F)^{-1}F^tC^{-1}z,
\end{equation}
where it can be seen that one matrix-vector and one matrix-matrix
products involving the inverse covariance matrix are required. When
$n$ is large, as the number of operations involved in the computation
of $C^{-1}$ scales with $n^3$, the computation of $Z^*$ is made
intractable, unless a particular form for $C$ is adopted. The
prediction variance or kriging variance can also be computed as:
\begin{align}
  \label{eq:krigingvar}
  \V(Z(x_0) - Z^*) = &C(0)-C_0^tC^{-1}C_0 +  \nonumber\\
&(f_0-F^tC^{-1}C_0)^t(F^tC^{-1}F)^{-1}(f_0-F^tC^{-1}C_0),
\end{align}
where $C(0) = \V(Z(x))$, corresponding to the punctual variance. We
remind that, \eqref{eq:uk} and \eqref{eq:krigingvar} are identical to the
conditional expectation and conditional variance under the Gaussian
hypothesis.

Two main families of approaches have been designed to overcome this
limit and are desribed in the following subsections. We then develop
an approach where they are combined.

\subsection{Covariance tapering}
Introduced in Furrer et al. \cite{tapering}, covariance tapering is an
approximation to the standard linear spatial predictor. The idea is to
taper the spatial covariance function to zero beyond a certain range
using a positive definite but compactly supported function. This terms
to deliberately introduce zeros into the matrix $C$ to make it
sparse. The linear system \eqref{eq:kriging} with a sparse (enough)
covariance matrix can be solved efficiently.

How the zeros are introduced is crucial. In particular, the positive
definiteness of any sparse modification of the covariance matrix must
be maintained. Let $\C_{\tau}$ be a covariance function that is
identically zero outside a particular range described by $\tau$. The
tapered covariance function is the product of $\C_{\tau}$ and the
initial covariance function $\C$:
\begin{equation}
\label{eq:tapering}
\C_{tap}(x,y)=\C(x,y)\C_{\tau}(x,y)
\end{equation}
The intuition behind this choice is both that the product $\C_{tap}$
preserves \emph{some} of the shape of $\C$ (in particular, its
behaviour near the origin) and that it is identically zero outside of
a fixed range. $\C_{tap}$ is a valid covariance function, as the
product of two positive definite functions. Finally, $\C_{tap}$ allows
to build a sparse covariance matrix that can be inverted using sparse
linear algebra techniques (e.g.  sparse Choleski decomposition), to
compute \eqref{eq:kriging} and \eqref{eq:krigingvar}.

The theoretical justification of tapering relies on a result for
kriging with misspecified covariance \citep{stein_equi}, stating
roughly that as far the behavior at the origin of the covariance is
well represented then the kriging predictor reaches optimality under
infill asymptotics settings, that is when the prediction is performed
on a bounded domain with increasingly denser data locations.

The main drawback of this approach is that it neglects the long range
structure that the data may present. While this has negligible effect
in densely sampled areas the effect can be significant in less sampled
areas. Further more, the range of the taper $\tau$ is set such that
the resulting covariance is sparse enough and needs to be reduced as
the number of data increases. Consequently, the prediction accuracy
can suffer from gaps between data locations when the sampling is not
uniformly dense.

\subsection{Low rank approaches}
Many low rank approaches have been proposed over the past few years,
among which we can cite the fixed rank approach \citep{fixedrank},
predictive processes \citep{pred_proc} and also more recent works
\citep{katzfuss2017multi,nychka2015multiresolution}. Precisely, under
this framework the covariance matrix of the observed random vector can
be written:
\begin{equation}
\label{eq:fixr}
C\,=\,\sig^2I\,+\,P\,B\,P^t,
\end{equation}
where $I$ is the identity, $P$ a known $n \f p$ matrix, $B$ a $p \f p$
covariance matrix. Such a model naturally arises by considering a
spatial process $Z$ such that
\[ Z(x) = \sum_{i=1}^p \eta_i\,P_i(x),\]
where the $P_i$ are known basis functions and
$(\eta_1,\ldots,\eta_p)^T$ is multivariate normal with mean 0 and
covariance $B$.

$p$ must be sufficiently small, so that the inverse of $C$ can be
computed. Indeed, $C^{-1}$, via the Woodbury-Morrison formula, takes
the following form:
\begin{equation}
C^{-1}=\sig^{-2}I-\sig^{-2}P\{\sig^2B^{-1} +P^tP\}^{-1}P^t,
\end{equation}
which only involves inverting the $p \f p$ matrices, ensuring a fast
computation of \eqref{eq:kriging} and \eqref{eq:krigingvar}, as far as $p$
is small enough.

The inference the parameters $\sig^2$ and $B$ can be conducted by a
method of moments \citep{fixedrank} or by maximum likelihood
estimation through an expectation-maximization (EM) algorithm
\citep{katzfuss2011spatio,braham2017fixed}.

An important point is the choice of the basis functions $P_i$. A wide
set of options is available: Cressie and Johannesson \citep{fixedrank}
propose using bi-square functions, the W-wavelet basis functions are
considered in \citep{wnychka}, Stein \citep{stein_kor} uses a
normalized Legendre polynomials basis on the sphere. In the proposed
combination approach described below, we show that the choice of the
basis functions can be data driven, using penalized least squares.

Finally, it is important to emphasize that this approach will neglect
the small scale fluctuations of the signal of interest. Indeed, the
fine scale fluctuations are represented by an unstructured noise.
This generally results in a loss of prediction accuracy and
oversmooth prediction maps.

\subsection{Combination}
\label{subsec:combi}
Knowing the advantages and drawbacks of theses two approaches, it
seems natural to try to combine them, so that the variability of the
random field can be captured at all scales. This has been proposed in
\citep{full_scale_cov} where the remainder of a predictive process
approximation is tapered. We suggest here to proceed in the opposite
way, first considering the small scale fluctuations before
investigating the modeling of the larger scales.

The predictor can be computed for large datasets as the covariance
matrix is of the same form as in \eqref{eq:fixr}, where the identity
is replaced by a sparse matrix.  More specifically, the form $C =
A+PBP^t$, still allows to use the Woodbury-Morrison matrix identity
that takes then the following form:
\begin{equation}
\label{eq:smw}
C^{-1} = A^{-1}-A^{-1}P(B^{-1}+P^tA^{-1}P)^{-1}P^tA^{-1}.
\end{equation}
Furthermore, if $A$ is sparse, $P$ is $n\f p$ and $B$ is $p\f p$, with
$p \ll n$, the numerical complexity of \eqref{eq:smw} scales with
$p^3\f n $. Indeed, \eqref{eq:smw} involves only the inversion of the
sparse matrix $A$ and the $p\f p$ matrices $B$ and
$(B^{-1}+P^tA^{-1}P)$.


It is interesting to notice that the constraints imposed on the
covariance model by its scalability with the number of data lead to
the definition of a particular form of covariance function where the
fluctuations of the random field are decomposed among several
components:
\begin{equation}
\label{eq:combinug}
    \C(x,y) = \A(x,y)+P(x)^tBP(y)+\sig_x\sig_y\delta_{xy},
\end{equation}
where $\A$ is a continuous, compactly supported covariance function
(it vanishes when $|x-y|$ is larger than a given, fixed range), $B$ is
a symmetric $p\f p$ positive definite covariance matrix, $P$ is a set
of basis functions and $\sig_x$ is a, possibly nonstationary, nugget
effect.

Specifically, the $\A$ covariance function is built using the
construction used in covariance tapering, namely as the product of any
covariance function with a compactly supported one (e.g. the spherical
model):
\begin{equation}
\label{eq:suppcompcov}
  \A(x,y) = \C_\tet(x,y)\C_\tau(x,y),
\end{equation}
where $C_\tet$ is any covariance function, possibly non stationary
\citep{paciorek2006spatial, fouedjio2016generalized,porcu2009quasi}
and $\C_\tau$ is a compactly supported correlation function with a
short fixed range \citep{compcov}. This construction induces the
corresponding covariance matrix $A$ to be sparse, with a density
controlled by the fixed scale parameter of $\C_\tau$. This
construction is the same as tapering \citep{tapering} but its use is
somehow different: rather than bulding a compactly supported
approximation of a non compactly supported covariance function, we
only use a compactly supported covariance function to model the small
scale fluctuation.

The first and last term of \eqref{eq:combinug} allow to model the
short scale fluctuations of the phenomenon under study which is
necessary to interpolate accurately the data in densely sampled
areas. Indeed, as will be seen in section \ref{sec:app}, the low rank
term only in \eqref{eq:combinug} fails to capture the fine scale
structure. Moreover, for most classes of covariance models that is
when the sampling becomes denser and denser in a fixed domain, the
quality of the kriging interpolator is essentially guided by the short
scale behaviour \citep{steinbook}.


The second term of \eqref{eq:combinug} allows to capture the
intermediate scale fluctuations, between those modeled by $\A$ and by
the drift. This representation is non stationary by construction. The
scale of the fluctuations modeled by this term is determined by the
choice of the basis functions constituting $P$. Therefore, they must
be chosen with care and not too numerous so as to keep the inversion
numerically attractive.

Finally, the nugget component of \eqref{eq:combinug} is classically
used to model a measurement error, for instance due to the internal
variability of the sensor used, or to model a non-continuous behaviour
of the regionalized variable, hence the term "nugget". We can allow
the value of the nugget to vary smoothly across the domain of interest
and filter it out at observation locations in prediction.

\subsection{Random field model}
Assuming the covariance model of \eqref{eq:combinug}, the random field model
can be rewritten as follows:
\begin{equation}
\label{eq:complete}
Z(x) = \mu(x) + \sum_{i=1}^{p}\eta_iP_i(x) + S(x) + \sig_x\eps(x),
\end{equation}
where $\eta = (\eta_1,\ldots,\eta_p)^t$ is a vector of centered,
square-integrable random variables with covariance matrix $B$,
$(P_i)_{i=1,\ldots,p}$ is a set of basis functions, $S$ is a centered
second order random field with covariance $\A$ \eqref{eq:suppcompcov},
and $\eps(x)$ is a centered, square-integrable, non structured random
function with variance 1. Moreover, to respect the model
\eqref{eq:combinug} $S$, $\eta$ and $\eps$ have to be uncorrelated.
$Z$ can hence be decomposed into four terms, each corresponding to a
different scale of fluctuations.

\subsection{Which basis functions?}
In low rank approaches, the basis functions are usually set in
advance. As said in the introduction, there is a wide set of options
for the basis functions to be used but the choice remains
arbitrary. We propose a data driven way to select the basis functions
among a vast family of candidates, possibly including anisotropies.

First note that equation \eqref{eq:complete} can be rewritten as
\begin{equation}
\label{eq:residual}
  Z(x) -\mu(x)= \sum_{i=1}^{p}\eta_iP_i(x)+Y(x),
\end{equation}
where $Y(x) = S(x) + \sig_x\eps(x)$ has only short spatial
correlations. In other words, we can express the detrended data as a
regression model with a random effect. We propose therefore to operate
a data driven selection of basis functions from a large set of
candidates, using variable selection methods currently used in
variable selection problems such as the LASSO \citep{lasso}. This
point will be further detailed in the next section.

Thinking about the Karhunen-Loève expansion \citep{loeve} makes this
idea natural. Indeed, any centered second order random field admits
the following representation:
\begin{equation}
  \label{eq:kl}
  \sum_{i\in I}\xi_i\lam_i \phi_i(x),
\end{equation}
where $I$ is at most countable, $(\xi_i)_{i\in I}$ are standard
uncorrelated random variables, $(\lam_i)_{ i\in I}$ is the set of the
of the covariance function of the field and the $(\phi_i)_{i\in I}$ is
the set of the associated eigenvectors. The rate of decay of the
eigenvalues is fast for covariances with smooth behavior, see
e.g. \citep{moi1,moi2}. Conversely, when considering a realization of
a Gaussian random field, most of the $\xi_i$ take values close to zero
in the representation \eqref{eq:kl}, as standard Gaussian random
variables. Consequently, only a small number of components are
"relevant" in the representation of the realization. Consequently,
modeling this scale of variation with a low number of carefully
selected basis functions seems to make sense. 

\section{Inference and prediction}
\label{sec:inf}
In this section, we will assume that $\mu(x)=0$ without loss of
generality. Indeed, the ordinary least square (OLS) estimate of the
trend, can be computed and subtracted from the data to get a centered
dataset. In the context of a large spatial dataset, the additional
variance due to the use of a non efficient estimator is generally
negligible.

The parameters to be inferred are the covariance matrix $B$, the
parameters of the covariance function $\A$ and the nugget $\sig_x$. We
describe here a step by step inference method. We also propose an
approach to select the most relevant basis functions $P$ among a large
dictionary of candidates.

We will first consider that the covariance $\sig^2\A_\tet(x,y) =
\sig^2_\phi\mathrm{R}(x,y) + \sig_x\sig_y\delta_{xy} =
\Cov(S(x)+\eps(x),S(y)+\eps(y))$, where $\tet$ is the set of small
scale covariance parameters, including $\sig^2_\phi$, the parameters
the correlation function $\mathrm{R}$ and $\sig_x$. In the stationary
case, $\sig_x = \sig_\eps$ is constant over the domain. In other words,
we will group the nugget and the short scale fluctuations term, as the
inference related to these two terms will be conducted conjointly. The
special case of an observation error modeled by a nugget effect will
be dealt with in a dedicated paragraph.

We consider a large vector of data $z = (z(x_1),\ldots,z(x_n))$. Under
Gaussian hypothesis, the full log-likelihood writes
\begin{align}
\label{eq:lik}
  L_z(\sig^2,\tet,B) = &-\frac{1}{2}\ln \det( \sig^2A_ \tet + P^t BP) \nonumber\\
&+ \frac{1}{2 \sig^2}z^t\left(A_ \tet^{-1}-\frac{1} \sig^2A_ \tet^{-1}P( B^{-1}+\frac{1} \sig^2P^tA_ \tet^{-1}P)^{-1}P^tA_ \tet^{-1}\right)z.
\end{align}
Maximizing directly the log-likelihood function is not straightforward
and cannot be computed analytically. Making appear the coefficients
$\eta$ of the basis functions as latent variables allows to decompose
the likelihood, indeed:
\begin{align*}
\eta &\sim \No(0, B) \\
Z|\eta &\sim \No(P^t\eta,\sig^2A_\tet).
\end{align*}
Consequently, an EM algorithm can be implemented to estimate the
parameters. Nevertheless, to do this, the basis functions need to be
known beforehand. One solution could be to add a penalty term in
\eqref{eq:lik} but this would only have add more complexity. To
simplify the parameter inference, we therefore rely on a step by step
algorithm where we treat the different scales separately, thus
lowering the computational burden. The several consecutive steps are
presented in table \ref{tab:steps} and detailed in the following
paragraphs.

\begin{table}
  \caption{Inference steps}
  \label{tab:steps}
  \begin{tabular}{cl}
    \hline \\ 1. & Fit the parameters of the small scale component
    (range and variance) by \\ & variogram fitting \\ \\ 2. & Select
    the basis functions by the LASSO \\ \\ 3. & Estimate $B$ and
    $\sig^2$ with fixed $P$ by ML \\ \\ 4. & Change the penalty in
    2. and update the selection\\ \\ 5. & Estimate $B$ with fixed $P$
    and $\sig^2$ by ML \\ \\ \hline
  \end{tabular}
\end{table}

We could have stopped at step 3. in table \ref{tab:steps} while
keeping a fixed $\sig^2$ while estimating $B$. Estimating both
$\sig^2$ and $B$ within the EM step aims at balancing the explained
variations between the small scale component and the basis functions,
while mitigating the possible parameter estimation errors caused by
variogram fitting. Indeed, the fitted variogram of the short scale
component try to explain the whole observed variability, at the
considered lags, without taking into account the additional basis
functions. Consequently, the penalty is generally set at a too large
value at first guess and the number of basis functions selected is too
low. Updating the value of $\sig^2$ within the EM algorithm allows
assigning a better weight to the short scale variations. A new
selection of basis functions can then be performed in accordance with
this new weight and the covariance matrix of its random coefficients
is then estimated. It is possible to iterate the process, updating
$\sig^2$ and redoing a selection, but this does not prove useful in
practice.

\subsection{Inference of the small scale structure}
The first step of the inference approach consists in fitting the,
possibly non stationary, small scale covariance structure $\A_\tet$.

When considering a stationary compactly supported covariance to model
the small scale structure, several options are available for the
estimation of the parameters. Maximum likelihood can be considered as
well as variogram fitting. The latter proves to be less computer
demanding and will therefore be used in the application section.

In presence of short scale non stationarity, a practical way to infer
non stationary covariance structures of the form given in
\citep{paciorek2006spatial} consists in fitting stationary models
locally as proposed in \citep{fouedjio2016generalized} or in
\citep{parker2016fused}. Interpolating the fitted value of the
parameters over the domain of interest allows to compute the value of
the covariance function between any two points of the domain, which is
necessary to compute the kriging interpolator. This approach only
requires that the parameters governing the non-stationarity vary
regularly over the domain, precisely the associated random function,
$S$ in our model, can be considered to be locally second order
stationary. This approach can however be time consuming, though an
embarrassingly parallel workload, and the benefits offered by
considering a non stationary small scale structure may not be worth
the trouble.

Recalling that the short scale structure should provides a sparse
matrix, the range of the compactly supported correlation function, the
taper, has to be small enough. We refer the interested reader to
\citep{tapering} for the choice of the taper. Basically, it has to be
chosen at least as smooth at the origin as the original model. 

Once the parameters of this structure have been inferred, the
remaining second order variability will be modeled using basis
functions. Indeed, the short scale structure explains a given part of
the variations of the phenomenon under study and this will guide the
choice of the basis functions to be used.
 
\subsection{Selection of the basis functions}
To estimate the matrix $B$, the basis functions constituting $P$ have
to be known. Numerous families of functions are available: local
bi-square functions \citep{fixedrank}, wavelets \citep{wnychka},
Legendre polynomials \citep{stein_kor}, smoothing splines, radial
basis functions, cosines, etc.

No clue exists however to choose the correct one for the problem at
hand, unless a covariance model has already been fitted and predictive
process can be considered \citep{pred_proc,katzfuss2017multi}.

To alleviate this choice, we propose to select the most relevant basis
functions among a large choice of candidates. This approach has
already been proposed in \citep{hsu2012group}, considering however
that the random variables $(\eta_i)_{i=1,\ldots,p}$ associated to each
basis function are independent, resulting in a diagonal $B$, which
seems too strong an assumption. The estimation of $B$, the covariance
matrix of the $(\eta_i)_{i=1,\ldots,p}$, is treated in the next
paragraph.

Here we propose as basis functions a collection of kernels anchored at
a set of knots homogeneously distributed over the domain. We choose to
use compactly supported covariance kernels for several reasons:
\begin{itemize}
\item we can build a multiresolution representation by considering
  several scale parameter values,
\item we can easily introduce anisotropic basis functions,
\item we can consider several levels of regularity at the
  origin (e.g. spherical, wendland, cubic, etc. kernels),
\item compact support makes the storage lighter and the prediction
  faster.
\end{itemize}

Practically, we consider a handful of values for the scale parameter
and orientation of the anisotropy. To each value of the scale
parameter corresponds a grid where the kernel is anchored. Scale
parameter values and orientations are crossed so as to generate
anisotropies. The scale parameter values are meant to correspond to
ranges greater than the one considered in the small scale structure,
so as to achieve scale separation between the two components.

Going back to \eqref{eq:residual}, the problem of selecting the relevant
basis functions among a dictionary of candidates terms to solve the
following penalized regression problem:
\begin{equation}
  \label{eq:select}
  \arg \min_{\be}\|(z-\widehat{\mu})-P\beta\|_2^2+\lam\|\beta\|_1,
\end{equation}
where $\be$ is a vector of estimates of the $(\eta_i)_{i=1,\ldots,p}$
and $\lam$ is a penalty term. The rationale for this is that if the
covariance of the random field admits a low rank representation, then
the random field can be expressed as a random linear combination of
the basis functions. Solving \eqref{eq:select} makes it possible to
select the active components in this representation. We choose to use
a $\mathrm{L}_1$ penalty so as to proceed to a more drastic selection
than that would have been obtained with the $\mathrm{L}_2$ norm
\citep{elements}.

Practically, \eqref{eq:select} is minimized though a coordinate descent
algorithm such as implemented in the \verb?R? package \verb?glmnet?
\citep{glmnet}. This type of algorithm has proven to be numerically
efficient for solving convex optimization problems such as
\eqref{eq:select}. The value of the penalty term is set by
cross-validation. As we already know, at this step, the variance
$\sig^2$ explained by the small scale structure, we propose to set the
penalty at the value necessary to explain the remaining
variance. Further more, as the regularization path is computed for a
vector of penalty values, we do not need to recompute the LASSO
when the value of $\sig^2$ is updated. This will be further explained
in section \ref{sec:app}.

It is worth noting that we are not interested in the estimated values
of $\be$, rather by the relevance of the associated basis function in
the explanation of the fluctuations of $z-\mu$. Furthermore, as the
remaining terms $S$ and $\eps$ deal with smaller scale structure, they
do not interfere in this selection procedure. For huge datasets, the
computational effort of this step can be alleviated using stochastic
gradient algorithms where the optimization is conducted iteratively
over subsets of data, see e.g.  \citet{bottou2010large}.

\subsection{Inference of the medium scale structure}
\label{subsec:em}
Once the $p$ basis functions have been selected, the covariance matrix
of their components has to be estimated and the weight of the small
scale structure updated. As proposed in \citep{katzfuss2011spatio}, we
choose to perform maximum likelihood estimation through an
expectation-maximization (EM) algorithm. The EM objective function
writes:
\begin{align}
\label{emfunc}
  Q(\Tet,\Tet_l) = &-\frac{n+p}{2}\ln(2\pi)
  -\frac{1}{2\sig^2}z^tA_\tet^{-1}z +
  \frac{1}{\sig^2}z^tA_\tet^{-1}P\mu-\frac{1}{2}\ln \det(\sig^2
  A_\tet)\nonumber\\ &
  -\frac{1}{2}\tr\left((\frac{1}{\sig^2}P^tA_\tet^{-1}P+B^{-1})(C+\mu\mu^t)\right)-\frac{1}{2}\ln\det(B),
\end{align}
where $\tr (X) $ is the trace of matrix $X$, $C = \V(\eta|Z=z,\Tet_l) =
\left(\frac{1}{\sig_l^2}P^tA_\tet^{-1}P+B_l^{-1}\right)^{-1}$ and $\mu
= \E(\eta|Z=z,\Tet_l) = \frac{1}{\sig_l^2}CP^tA_\tet^{-1}z$.

The update formulas for the parameters $\sig^2$ and $B$ are:
\begin{align}
     \sig_{\new}^2 &= \frac{1}{n}\E\left((Z-P\xi)^tA_\tet^{-1}(Z-P\xi)|Z,\Tet_l\right)\nonumber\\ 
&=\frac{1}{n}\left(Z^tA_\tet^{-1}Z -2Z^tA_\tet^{-1}P\mu + \tr((P^tA_\tet^{-1}P+B^{-1})(C+\mu\mu^t))\right),\\
    B_{\new} &= C+\mu\mu^t.
\end{align}

The algorithm is initialized with the value of $\sig^2$ obtained at
the first step and $B$ as the identity and run until a stopping
criterion is met. The stopping criterion we use is defined as:
\begin{equation}
  \label{eq:stopping}
  \min\left(|\sig^2_{l+1}-\sig^2_{l}|,\frac{|\sig^2_{l+1}-\sig^2_{l}|}{\sig^2_{l+1}}\right) + \min\left(\left|\log\left(\frac{\det(B_{l+1})}{\det(B_{l})}\right)\right|,\left|\frac{\log\left(\frac{\det(B_{l+1})}{\det(B_{l})}\right)}{\log\det (B_{l+1})}\right|\right).
\end{equation}
We decide to stop the algorithm when the criterion is less than a
given tolerance value for 20 iterations in a row.

\subsection{Update of the selection and final model}
In the previous step, the amount of variability explained by the small
scale structure modeled by $\sig^2$ is updated. It can be lowered or
increased depending on the data. It is thus possible, and sometimes
desirable, to update the selection of the basis functions accordingly.

As said earlier, the whole regularization path is computed in the
selection step. Consequently, when the value $\sig^2$ is updated in
the EM algorithm, there is no need for any further computation to
update the selection. The penalty value has only to be set so that the
basis functions explains the updated remaining variance.

Then, the EM algorithm is run once again to get the final estimate of
the matrix $B$, while the value of $\sig^2$ is left unchanged. It is
possible to iterate the process, updating $\sig^2$ and redoing a
selection, but this does not prove useful in practice.

\subsection{Prediction}

It is worth noting that most of the quantities needed for the
prediction have already been computed in the EM step. Indeed, under
the Gaussian hypothesis, the kriging predictor is identical to the
conditional expectation and the kriging variance to the conditional
variance. These quantities can be straightforwardly computed knowing
$\P(\eta|z)$. Consequently the kriging predictor at a new location
$x_0$ writes:
\begin{equation}
  \label{eq:condexp}
  \E(Z(x_0)|z) = P_0^t\mu + A_0A^{-1}(z-P^t\mu),
\end{equation}
where $S_0$ is the small scale component value at $x_0$, $\mu =
\E(\eta|z)$, $P_0$ is the vector of the basis functions at $x_0$ and
$A_0$ is the covariance of the small scale structure between the
target and the observation. Similarly, the kriging variance writes:
\begin{equation}
  \label{eq:condvar}
  \V(Z(x_0)|z) = \sig^2 - A_0A^{-1}A_0 - 2A_0A^{-1}P^tCP_0 + P_0^tCP_0,
\end{equation}
where $C = \V(\eta|z)$.  These formulas provides the most efficient
way to compute these quantities. They only require to compute the
cholesky decompositions of the sparse matrix $A$ and the small matrix
$C$. It is interesting to note that although the covariance model is
built so that we can use the Woodbury-Morrison formula to solve the
kriging system, this formula is not used explicitly to compute the
predictor and the prediction variance.

\subsection{Filtering a measurement error}

When an observation error exists, it is important to filter it out
from the prediction, including at observation locations. In that case,
the observation error is modeled by a nugget effect. A first idea
could be to fit the variance of the nugget at the variogram fitting
step and to keep it fixed in the update of the variance term in the EM
step. In other words, we would have $\sig^2\A_\tet(x,y) =
\sig^2_\phi\mathrm{R}(x,y) + \sig^2_{\eps}\delta_{xy}$, where
$\sig_{\eps}$ is fixed. When doing this unfortunately, there is no
closed form expression for the update formula of $\sig^2$. Therefore a
numerical optimization is needed at each iteration of the EM, which
slows down the calculations.

Here we propose a alternative scheme. We first rewrite the small scale
covariance as $\sig^2\A_\tet(x,y) = \sig^2_\phi\left(\mathrm{R}(x,y) +
\frac{\sig^2_{\eps}}{\sig^2_\phi}\delta_{xy}\right)$. Then, we fix the
value of $\frac{\sig^2_{\eps}}{\sig^2_\phi}$ as its first estimates
$\frac{\widehat{\sig^2_{\eps}}}{\widehat{\sig^2_\phi}^{\mathrm{old}}}$
prior to go through the inference steps described above, where we get
a new estimate $\widehat{\sig^2_\phi}^{\mathrm{new}}$. Finally, as the
nugget variance is generally correctly fitted at the variogram fitting
step and we do not want to change it, we reinject it in the estimated
small scale covariance so that it takes the following form:
\begin{equation}
  \widehat{\sig^2}\A_{\widehat{\tet}}(x,y) =
  \widehat{\sig^2_\phi}^{\mathrm{new}}\left(1+\frac{1}{\widehat{\sig^2_\phi}^{\mathrm{old}}}-\frac{1}{\widehat{\sig^2_\phi}^{\mathrm{new}}}\right)\mathrm{R}(x,y)
  + \widehat{\sig^2_{\eps}}\delta_{xy}.
\end{equation}
By doing this, we keep constant the nugget effect and update the
weight of the small scale component. We can reasonably assume that
this simplification has a negligible effect on the computation of the
likelihood in the EM algorithm.

Finally, to filter the error from the observation, the kriging is
performed considering the following small scale covariance matrix:
\begin{equation}
  \widehat{\sig^2_\phi}^{\mathrm{new}}\left(1+\frac{1}{\widehat{\sig^2_\phi}^{\mathrm{old}}}-\frac{1}{\widehat{\sig^2_\phi}^{\mathrm{new}}}\right)\mathrm{R}(x,y).
\end{equation}

\section{Application}
\label{sec:app}
In this section, we provide experiments of the proposed method and
comparison with existing approaches, on two simulated and one real
datasets. The different approaches have been implemented in \verb?R?
with an extensive use of the packages \verb?spam? \citep{spam},
\verb?glmnet? \citep{glmnet} and \verb?RGeostats?
\citep{rgeostats}. The computations have been performed a laptop
computer with 3rd generation Intel Core i7-3687U processor (year 2013)
and 8 GB of memory.

\subsection{Nested covariance}
In this first example, we consider a simulated Gaussian
random field (GRF) on $[0,1]^2$, centered with covariance
\begin{equation}
  \C(x,y) = 0.5\text{sph}(|x-y|,0.05)+0.5\exp(|x-y|,0.1),
\end{equation}
where sph stands for the spherical covariance function, $\exp$ for the
exponential covariance function and the arguments are the distance and
the scale parameter. This kind of nested covariance is often used in
practice to model variables that exhibit variations at two different
scales

\begin{figure}[h!]
  \centering
  \begin{subfigure}[t]{0.45\linewidth}
    \includegraphics[width=\textwidth]{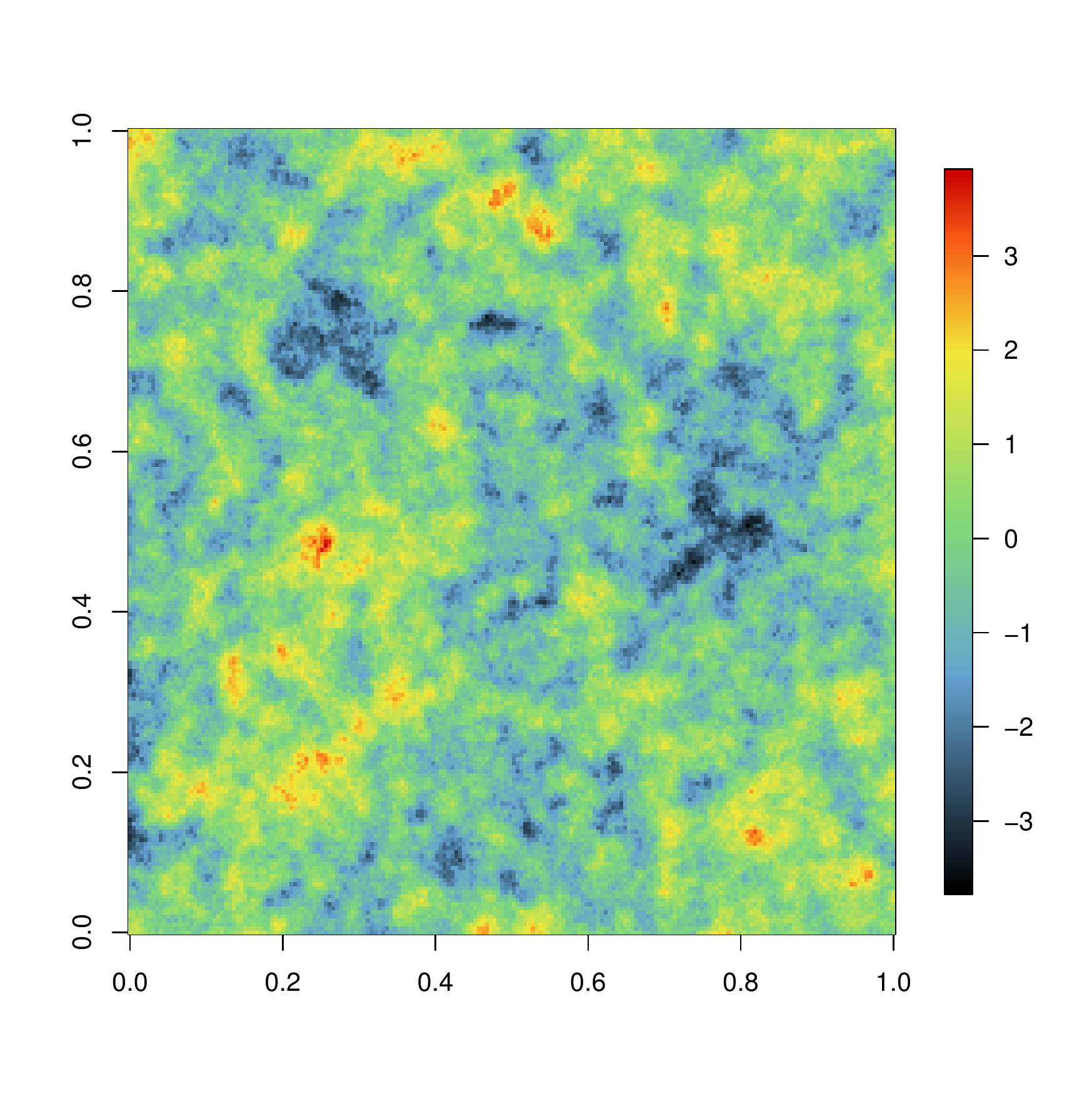}
    \caption{Reference realization}
    \label{fig:simgigref}
  \end{subfigure}
  \begin{subfigure}[t]{0.45\linewidth}
    \includegraphics[width=\textwidth]{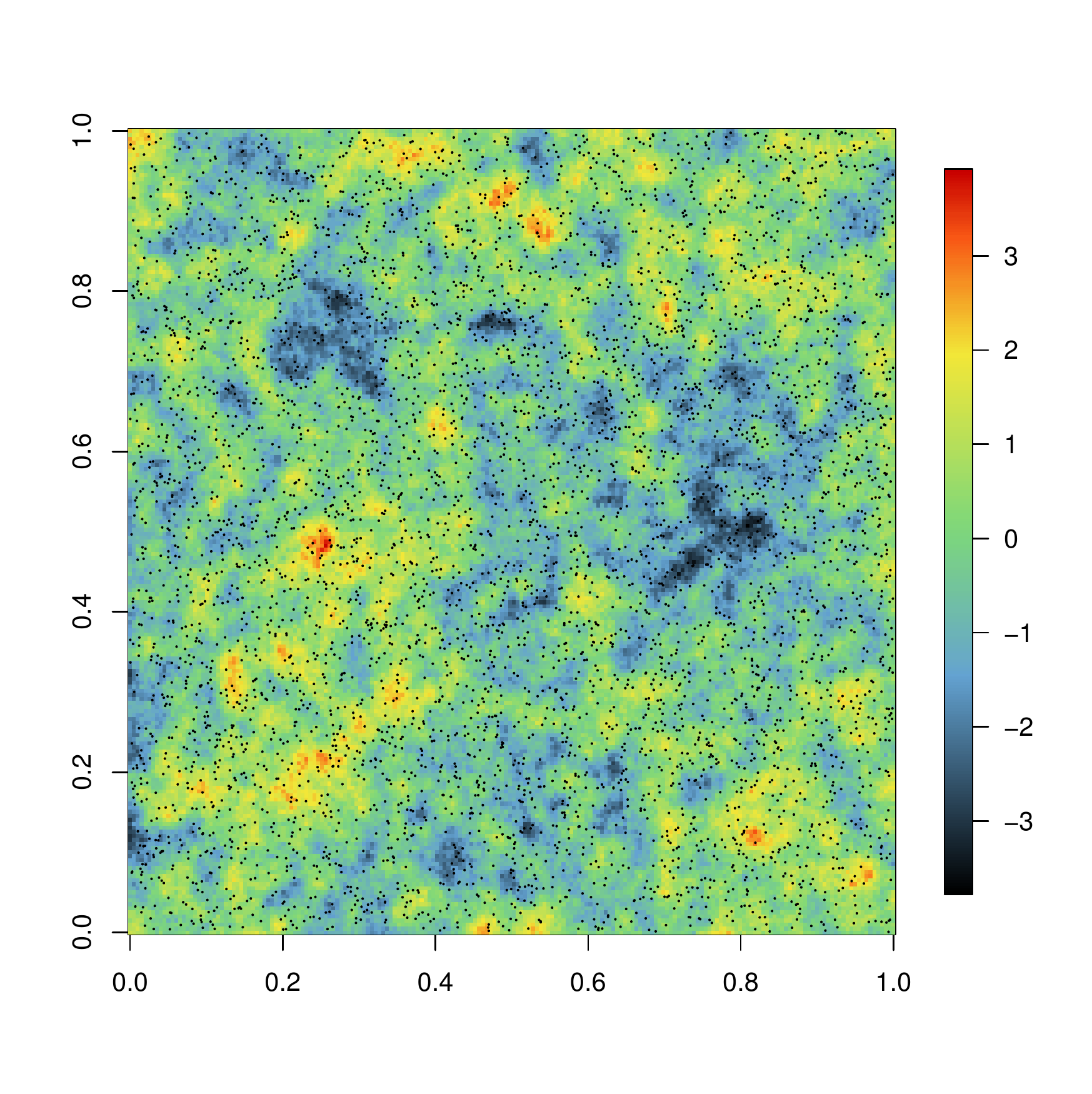}
    \caption{Sampling at 5000 locations}
    \label{fig:simgigsam}
  \end{subfigure}
\caption{Reference simulation and sampling}
\label{fig:simgig}
\end{figure}

The GRF is simulated on a 200x200 grid and sampled at 5000 points, as
shown in figures \ref{fig:simgigref} and \ref{fig:simgigsam}

\begin{figure}[h!]
  \centering
  \begin{subfigure}[t]{0.45\linewidth}
    \includegraphics[width=\textwidth]{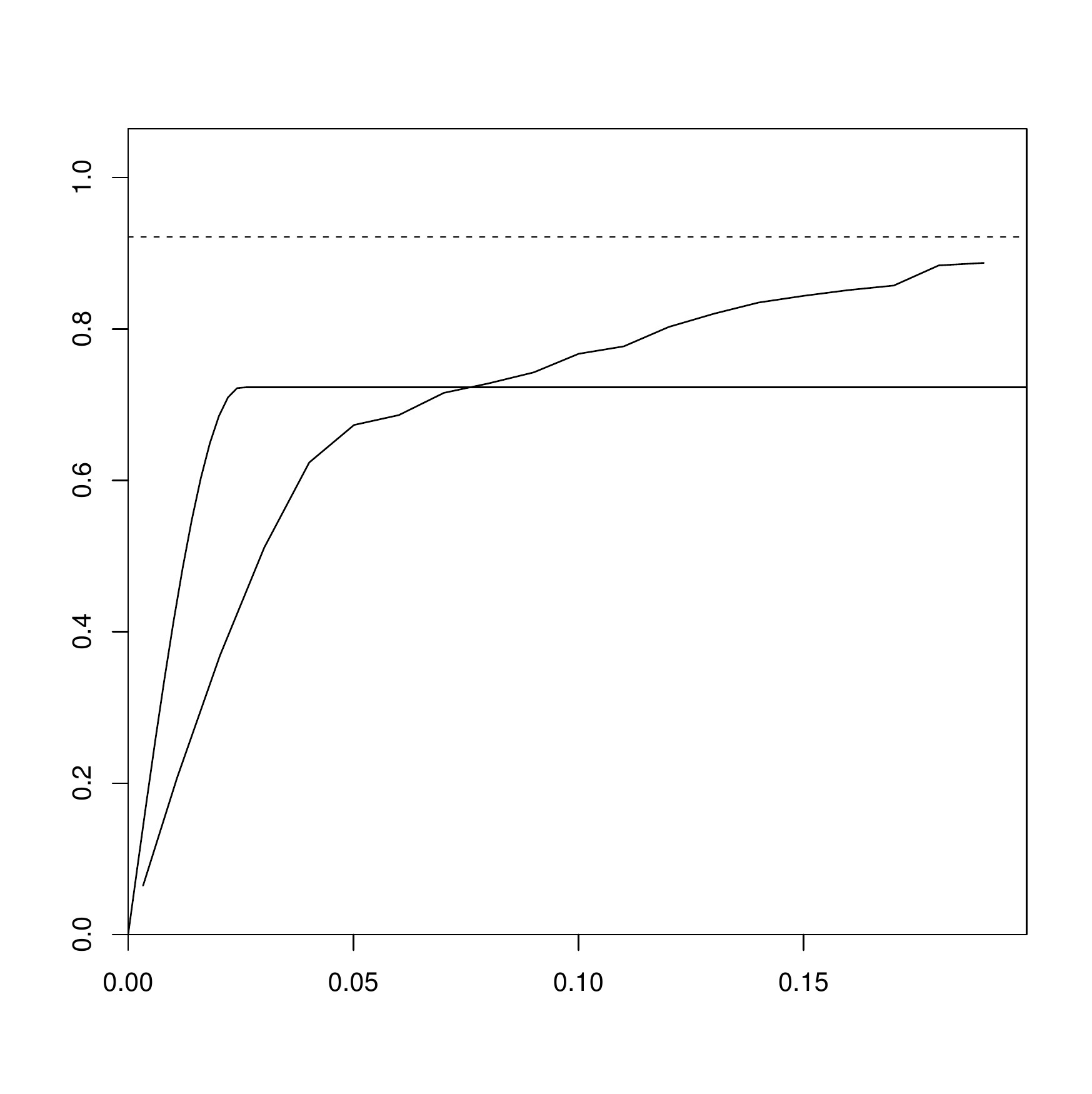}
    \caption{First fitting of the small scale structure}
    \label{fig:variogig1}
  \end{subfigure}
  \begin{subfigure}[t]{0.45\linewidth}
    \includegraphics[width=\textwidth]{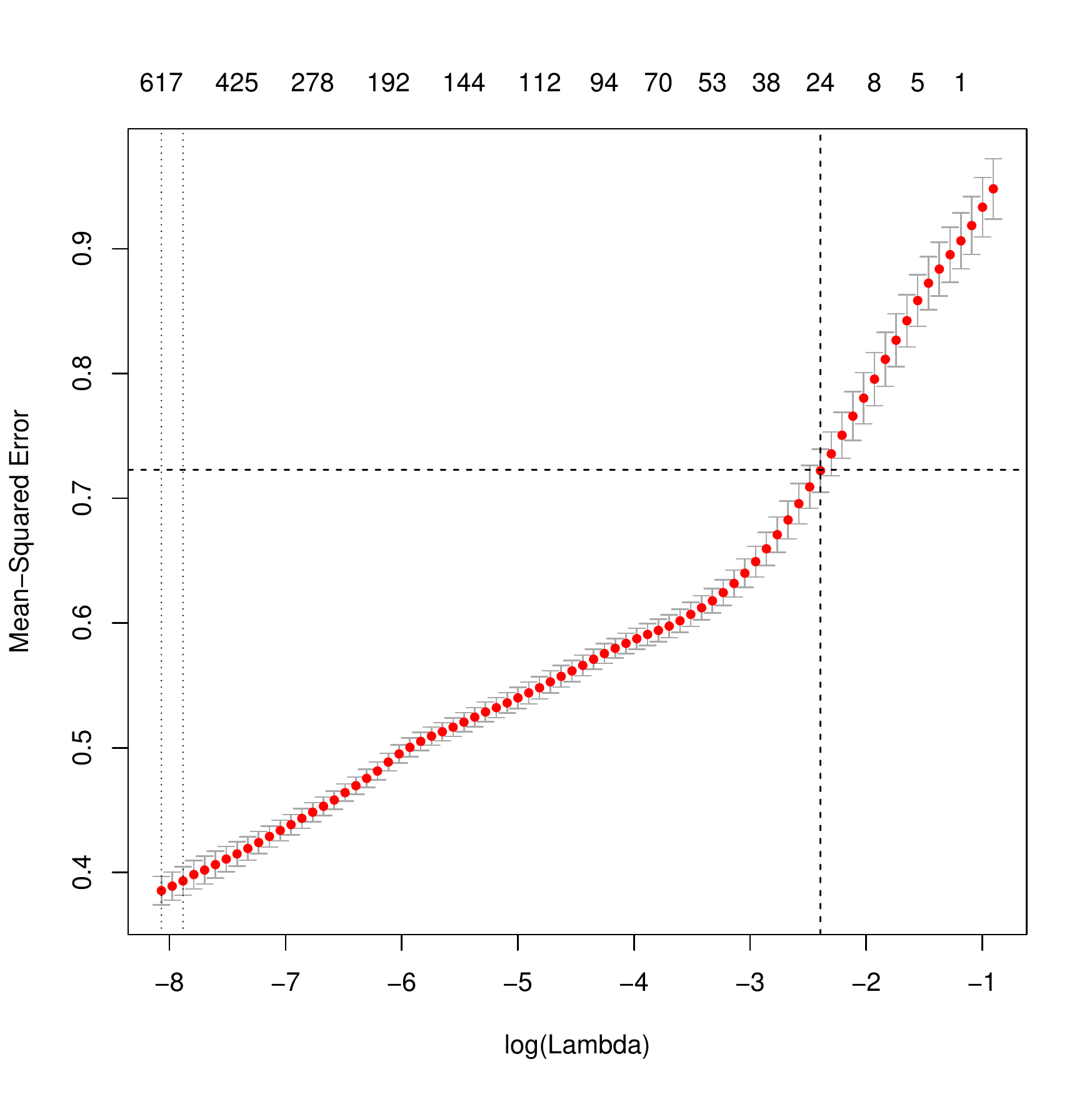}
    \caption{First selection of basis functions}
    \label{fig:glmgig1}
  \end{subfigure}
  \caption{First fitting of the model}
\end{figure}

From this dataset, we infer the combination model described in
paragraph \ref{subsec:combi} following the steps detailed in
\ref{sec:inf}. The first step consits in computing the empirical
variogram at short distances and fitting a tapered variogram with a
small range. This step is illustrated in figure \ref{fig:variogig1}.  The
structured inferred for the small scale component is:
\begin{equation}
  \label{eq:tapgig}
  A(x,y) = \sig^2\exp_{\tet}(x,y)\mathrm{sph}_{0.025}(x,y)
\end{equation}

The quantity $\sig^2$ can be interpreted as the amount of variability
explained by the short scale structure. As can be seen, from figure
\ref{fig:variogig1}, this value is overestimated but will be updated
later. To explain the remaining variability, we will select a set of
basis functions from a large set of candidates. The dictionary of
candidates is built combining several ingredients: two kinds of
compactly supported basis functions built from the cubic and spherical
covariance models, two ranges are considered for each one, 0.5 and
0.2, and four different anisotropy angles 0, $\pi$/4, $\pi$/2,
$3\pi$/4. The basis functions are computed over knots forming two
different grids whose discretization step is set according to the
range of the basis function.

This makes a total of 2658 candidates basis functions. The relevant
functions are selected by the LASSO using the \verb?cv.glmnet?
function of the package \verb?glmnet? \citep{glmnet}. Using this
function, a 5-fold cross-validation is performed so as to assess the
stability of the selection procedure. Then, the penalty is set so that
the basis functions explain the remaining variance, as shown in figure
\ref{fig:glmgig1}. We can see that 24 basis functions are selected as
a first guess. 

\begin{figure}[h!]
  \centering
  \begin{subfigure}[t]{0.45\linewidth}
    \includegraphics[width=\textwidth]{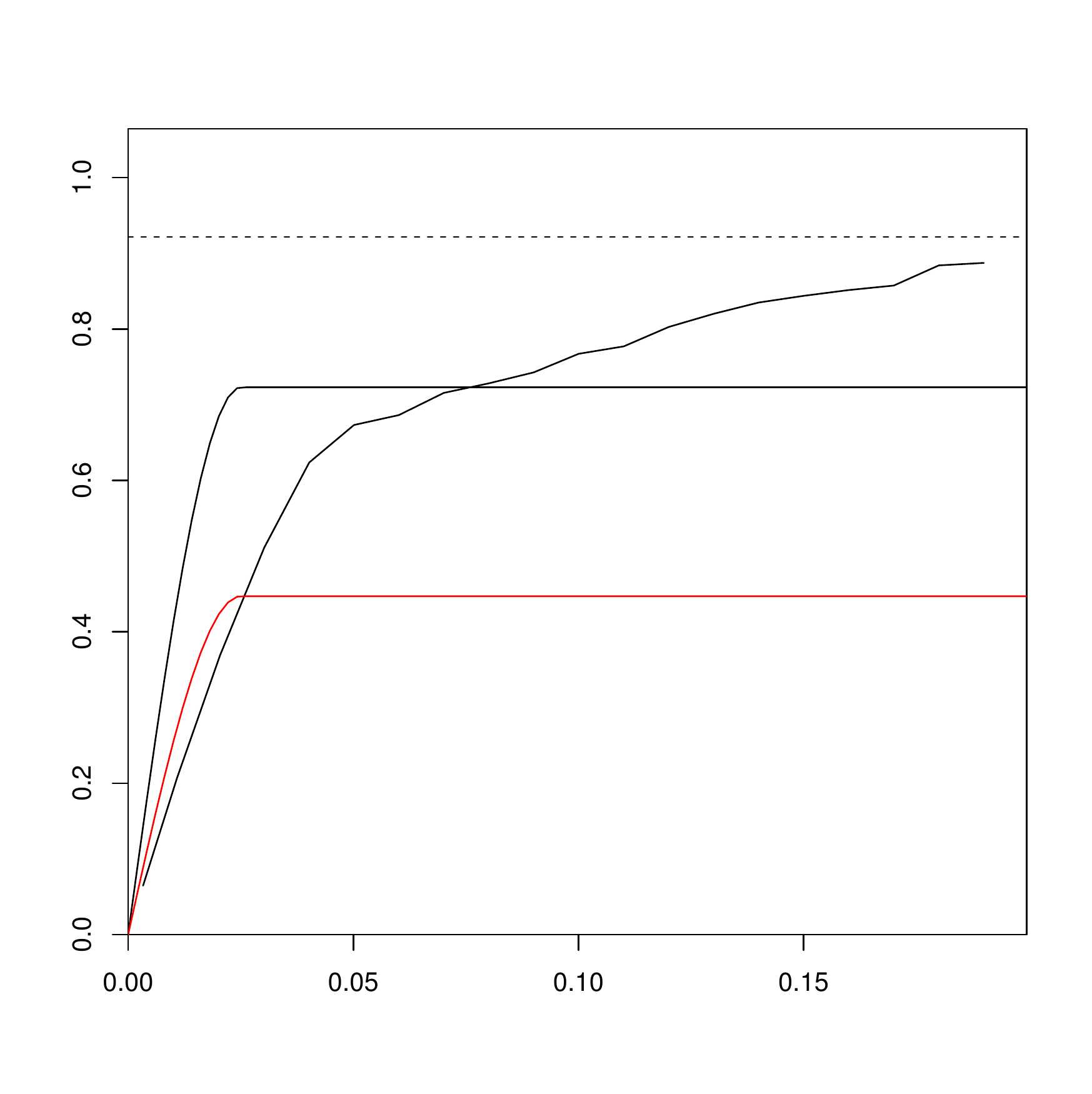}
    \caption{Updated variogram}
    \label{fig:variogig2}
  \end{subfigure}
  \begin{subfigure}[t]{0.45\linewidth}
    \includegraphics[width=\textwidth]{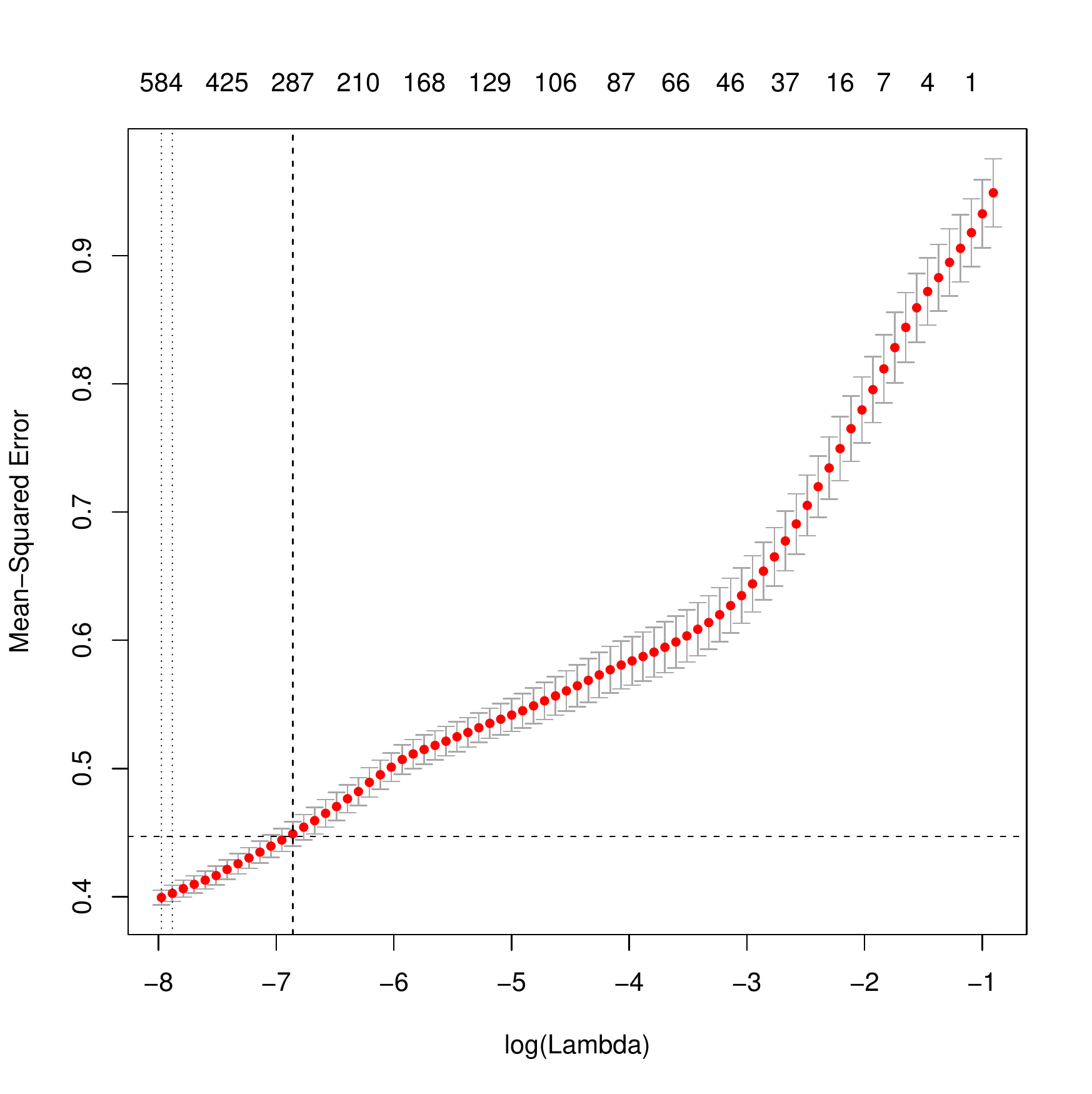}
    \caption{Updated selection of basis functions}
    \label{fig:glmgig2}
  \end{subfigure}
  \caption{Model update}
\end{figure}

The next step consists in running the EM algorithm to estimate the
covariance matrix $B$ of the random coefficients of the basis
functions and update the value of $\sig^2$. This allows to adjust the
part of the variations explained by the small scale component. The EM
algorithm is run for this selection until the stopping criterion is
met, with a tolerance value set at 10$^-3$. The updated variogram is
shown in figure \ref{fig:variogig2}, where we can see that the fit is
now more precise at small distances. A new set of basis functions is
selected accordingly, that represents 287 basis functions, as can be
seen in figure \ref{fig:glmgig2}.

\begin{figure}[h!]
  \centering
  \begin{subfigure}[t]{0.45\linewidth}
    \includegraphics[width=\textwidth]{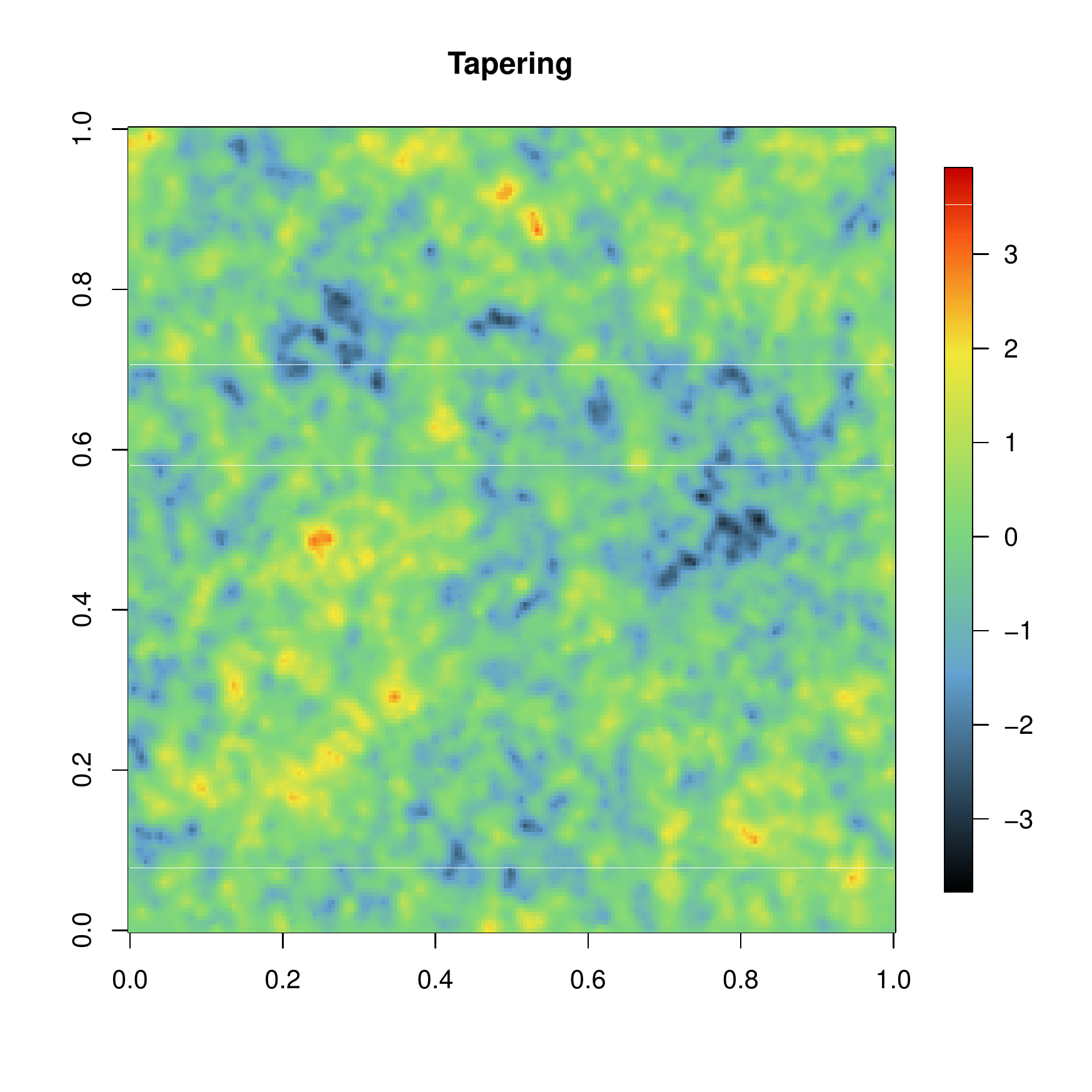}
    \caption{Covariance tapering}
    \label{fig:kriggig_tap}
  \end{subfigure}
  \begin{subfigure}[t]{0.45\linewidth}
    \includegraphics[width=\textwidth]{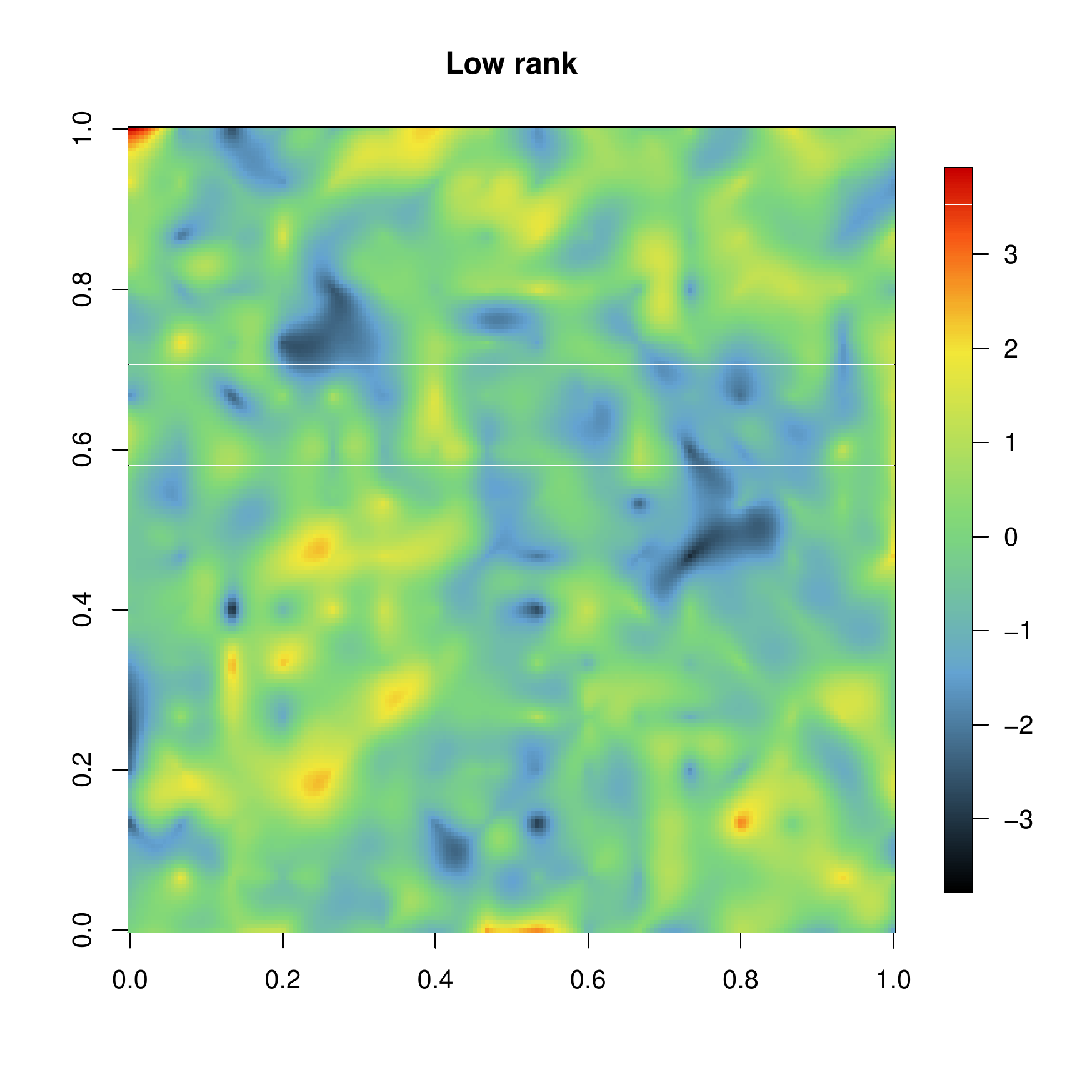}
    \caption{Low rank model}
    \label{fig:kriggig_lr}
  \end{subfigure}
  \begin{subfigure}[t]{0.45\linewidth}
    \includegraphics[width=\textwidth]{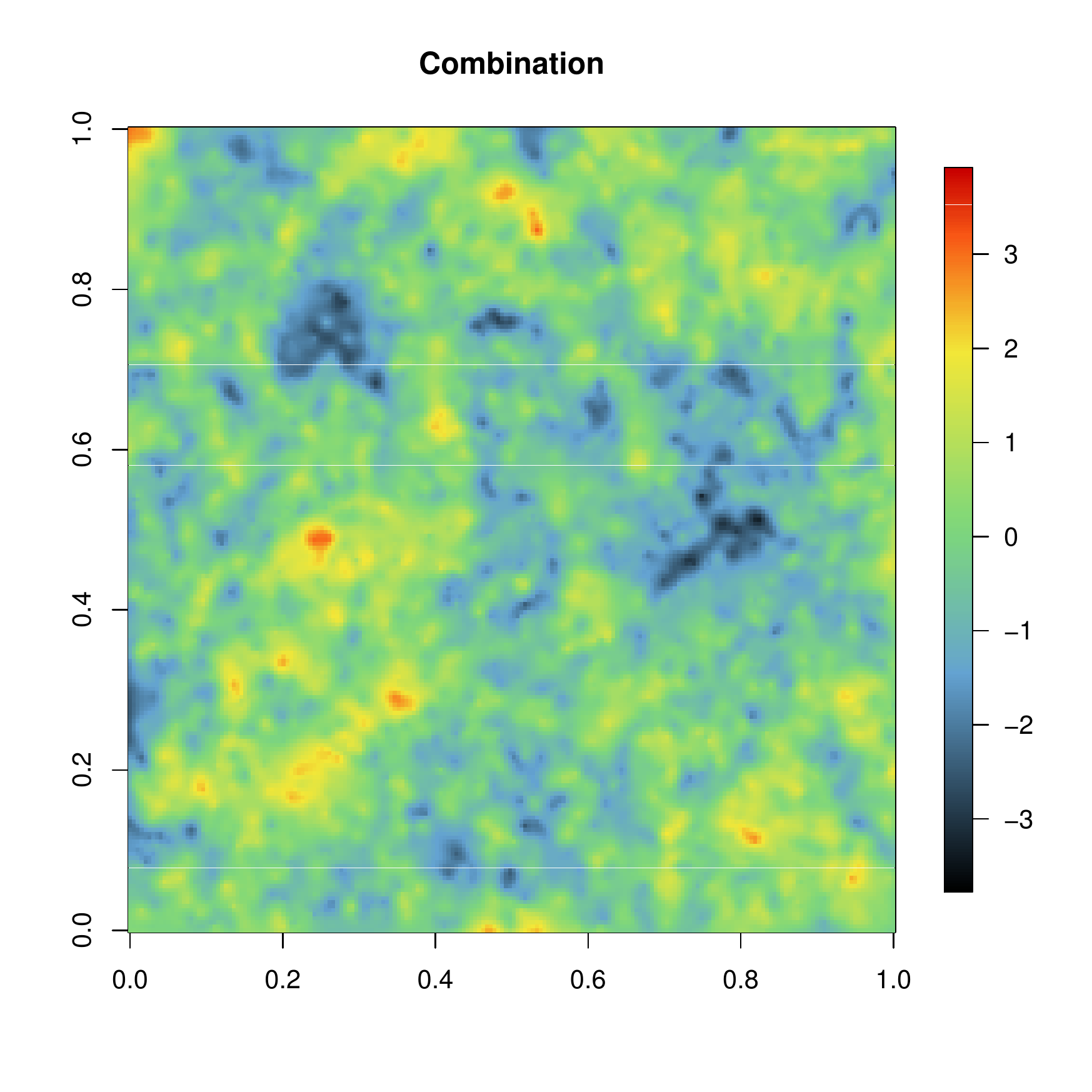}
    \caption{Combination}
    \label{fig:kriggig_combi}
  \end{subfigure}
  \begin{subfigure}[t]{0.45\linewidth}
    \includegraphics[width=\textwidth]{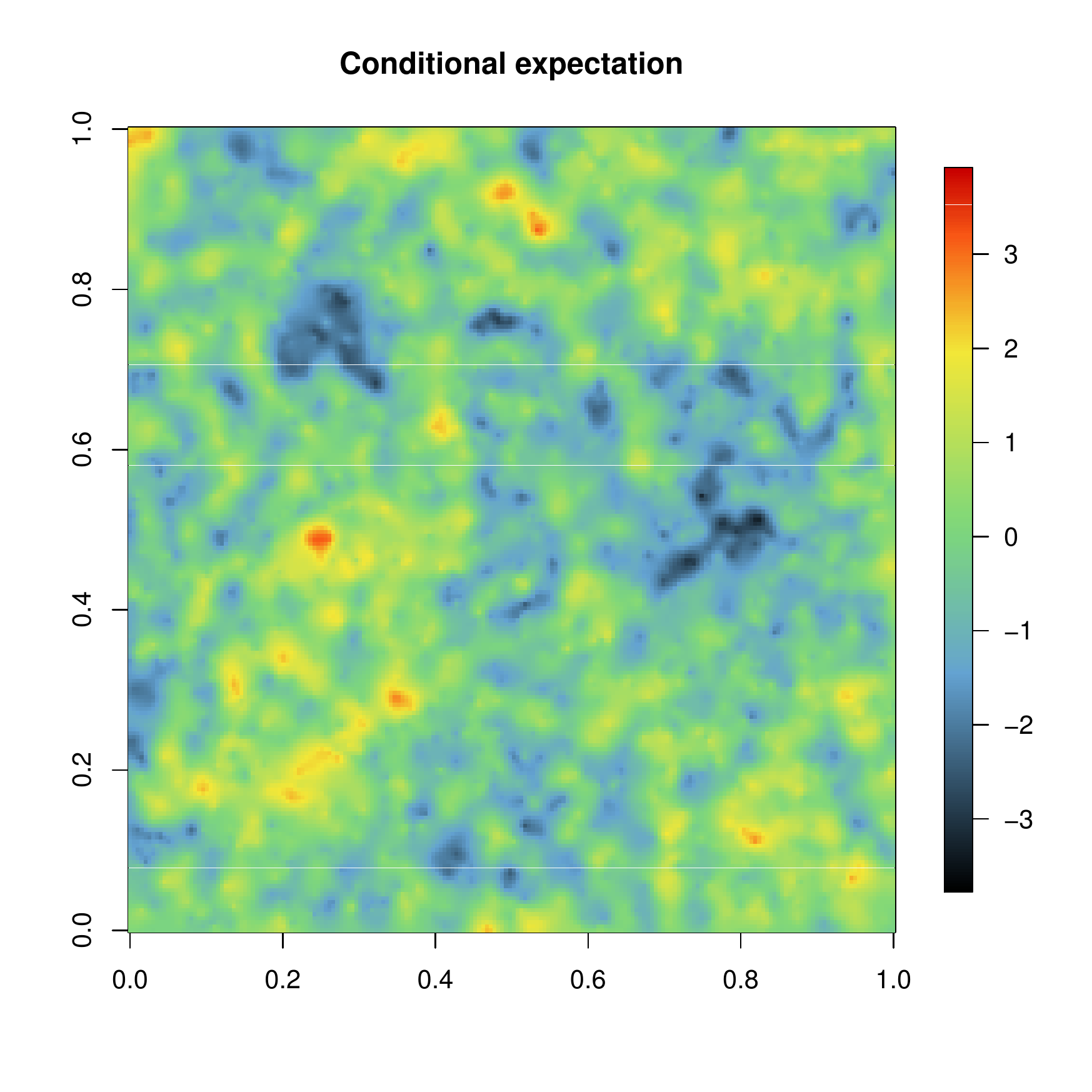}
    \caption{Conditional expectation}
    \label{fig:kriggig_ce}
  \end{subfigure}
\caption{Kriging maps}
\label{fig:kriggig_maps}
\end{figure}

The EM algorithm is run once again with the new set of basis functions
to estimate $B$, while the value of $\sig^2$ is kept fixed. Then the
kriging predictor can be computed on the 200x200 discretization grid
of the unit square. We also computed the map obtained using the
tapering term only \eqref{eq:tapgig}, the basis functions only, where
these are selected setting the penalty value at the largest value such
that the error is within 1 standard error of the minimum (560), and
the conditional expectation obtained using the true covariance. The
kriging maps are depicted in figure \ref{fig:kriggig_maps}.

While there are few visible differences between the maps obtained with
covariance tapering, the combination and the true model, it can be
seen in figure \ref{fig:kriggig_lr} that the map obtained with the low
rank model is much smoother. This is due to the bad representation of
the small scale fluctuations using only the basis functions. In that
case, the small scale fluctuations are modeled by an unstructured
noise whose variance is close to the $\sig^2$ term estimated in the
combination approach.

\begin{figure}[h!]
  \centering
  \begin{subfigure}[t]{0.3\linewidth}
    \includegraphics[width=\textwidth]{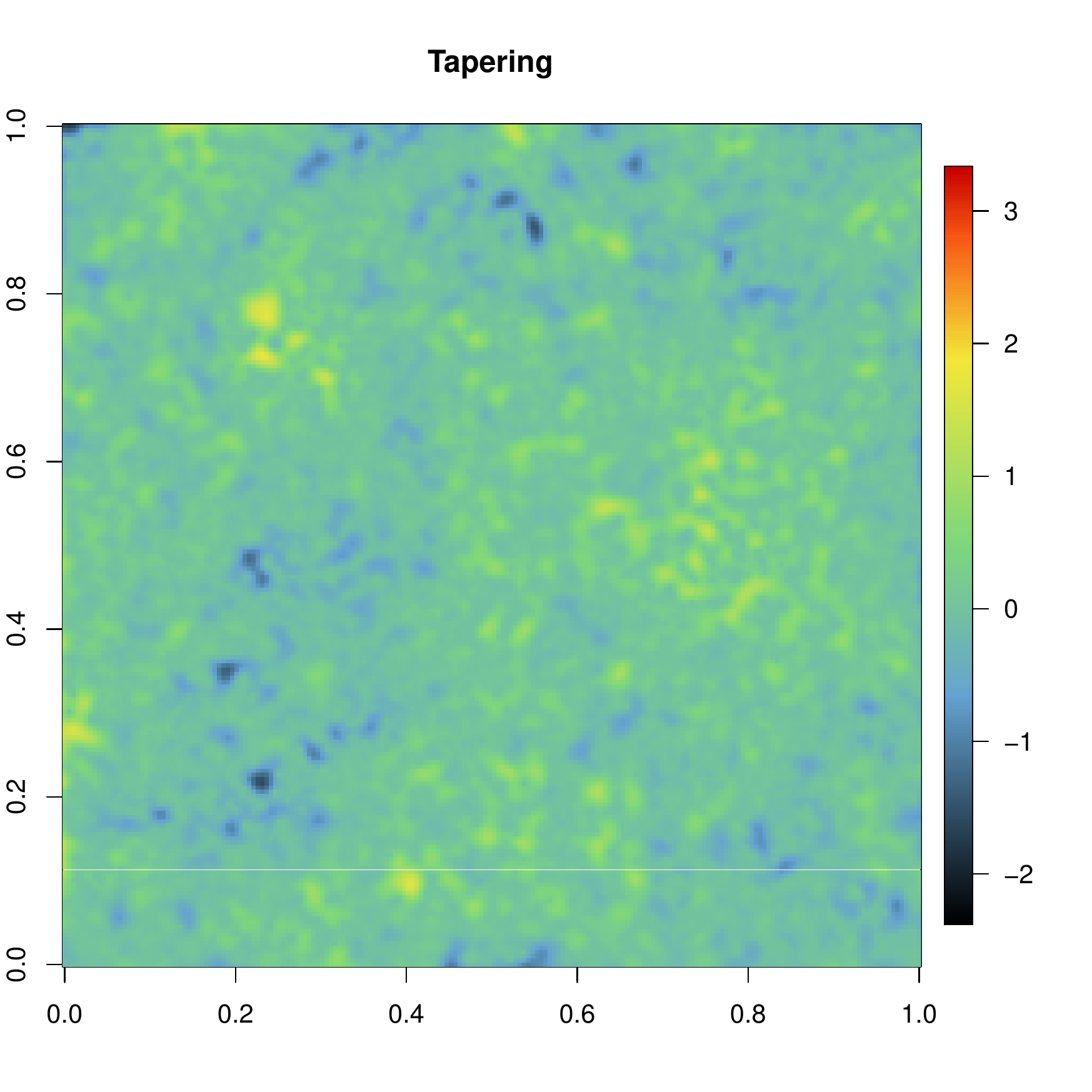}
    \caption{Covariance tapering}
    \label{fig:kriggigdiff_tap}
  \end{subfigure}
  \begin{subfigure}[t]{0.3\linewidth}
    \includegraphics[width=\textwidth]{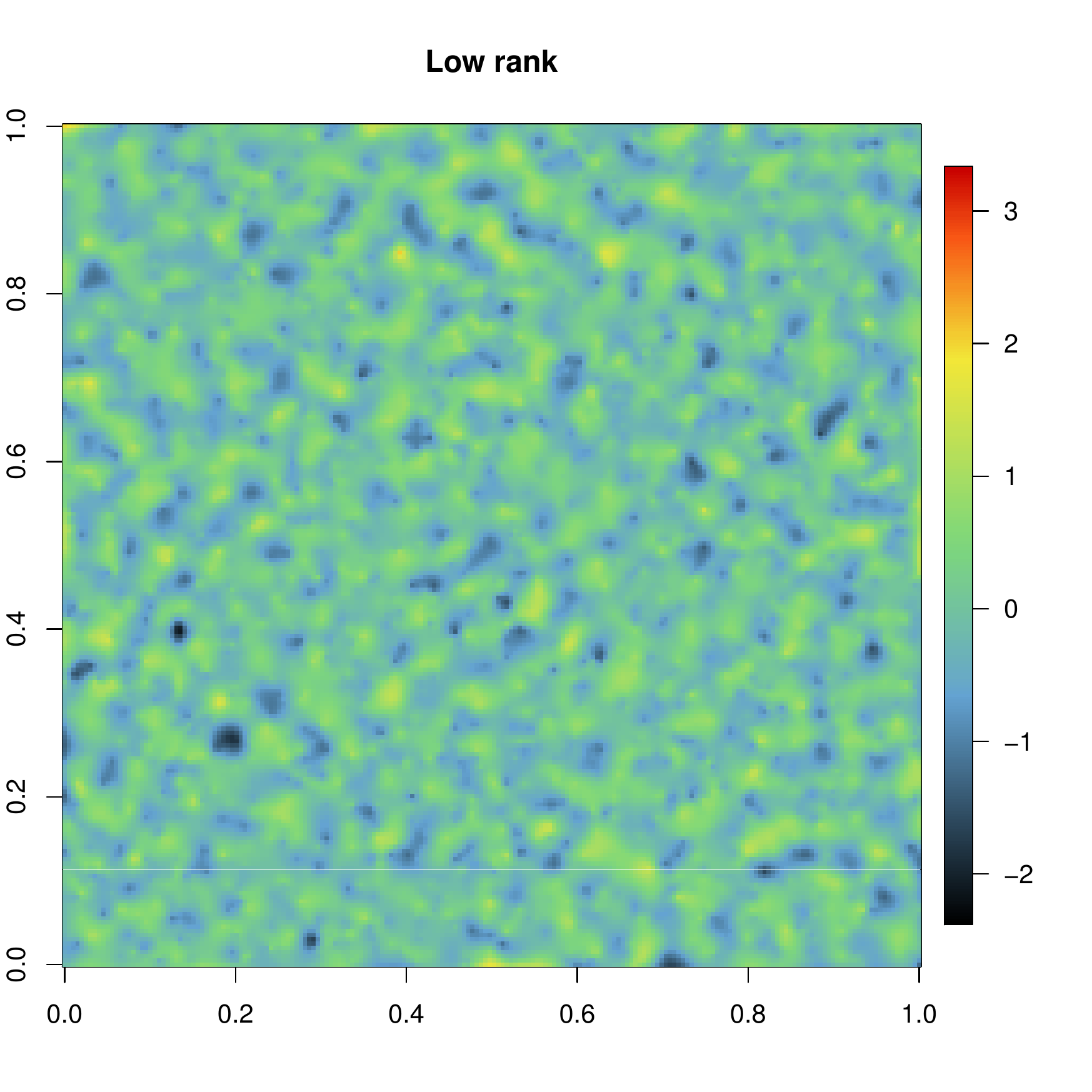}
    \caption{Low rank model}
    \label{fig:kriggigdiff_lr}
  \end{subfigure}
  \begin{subfigure}[t]{0.3\linewidth}
    \includegraphics[width=\textwidth]{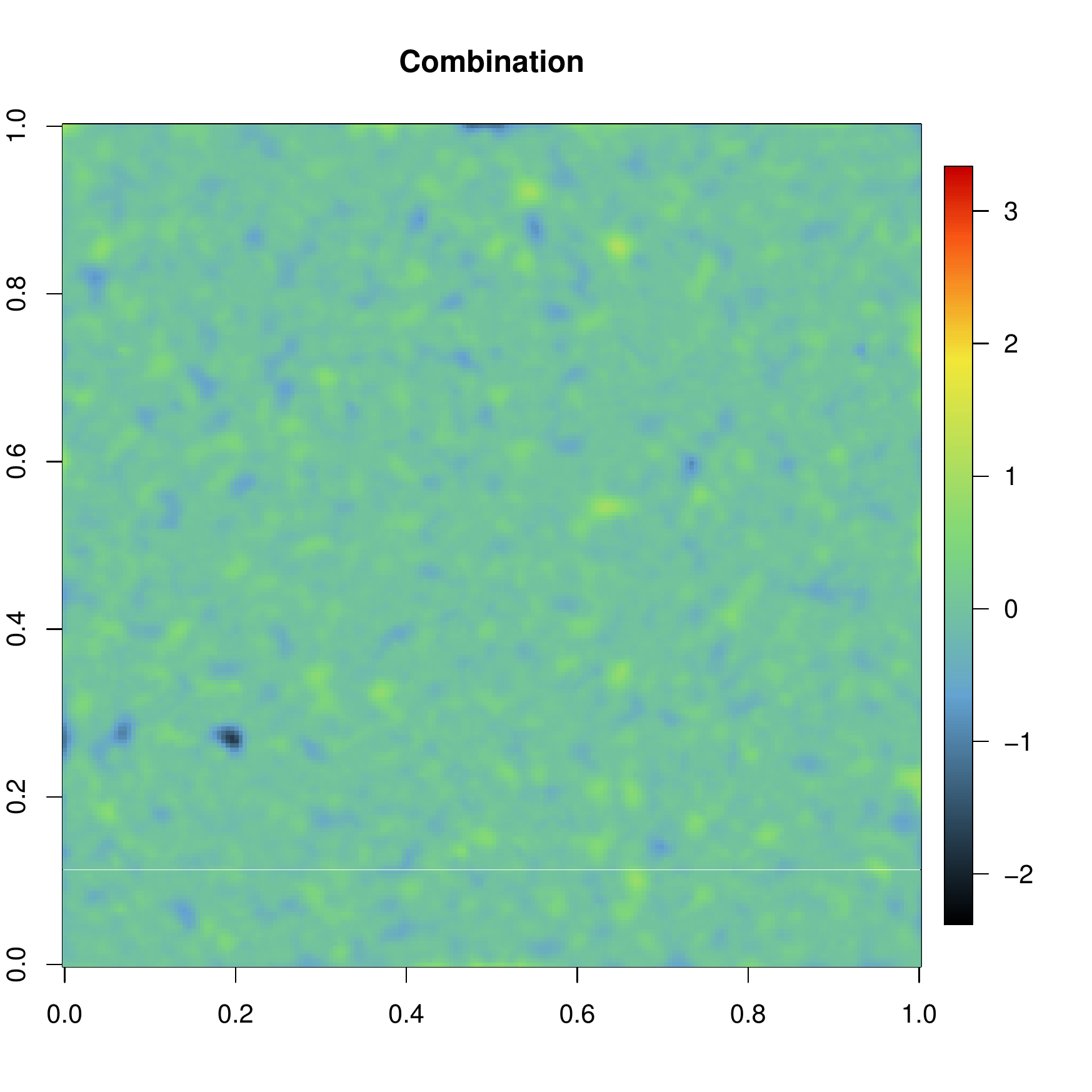}
    \caption{Combination}
    \label{fig:kriggigdiff_combi}
  \end{subfigure}
\caption{Difference maps}
\label{fig:kriggig_diff}
\end{figure}

We also computed the difference maps between the three fast approaches
and the conditional expectation (figure \ref{fig:kriggig_diff}).

Finally, we compute the mean squared prediction error (MSPE) obtained
by the four methods. They are shown in table \ref{tab:mspegig},
together with the details of each model. We can notice here that the
low rank approach obtains the worse result among the different methods
implemented, while using a large number of basis functions. As
explained earlier, it fails at capturing the small scale fluctuations,
modelling them as a nugget effect with a large variance and therefore
filters out this information while predicting, resulting in an
oversmooth kriging map.

\begin{table}
  \caption{Summary}
  \begin{tabular}{lccc}
    \hline\\
    model & small scale structure & number of basis functions & MSPE \\
    \hline\\

    Covariance tapering & range 0.025 & 0 & 0.24 \\

    Low rank  & nugget 0.40 & 560 & 0.35 \\

    Combination & range 0.025 & 287 & 0.19\\

    Conditional expectation & & &0.17\\
    \hline
  \end{tabular}
  \label{tab:mspegig}
\end{table}

\subsection{Non stationary Matérn covariance}
Recently, new non stationary covariance constructions have been
proposed, see
e.g. \citep{fouedjio2016generalized,porcu2009quasi}. These models
allow to make vary the anisotropy or even the regularity of the random
field over the domain considered. The inference of such model relies
on local methods (local variogram fitting or local likelihoods)
combined with an interpolation of the locally fitted parameters. This
procedure turns out to be computationally intensive. The resulting
dense covariance matrices are also cumbersome to compute. For our
second example, we consider a GRF on $[0,1]^2$, centered, with
covariance
  \begin{equation}
    \label{eq:covmatns}
    \Cov(Z(x),Z(y)) = \phi_{xy}\frac{2^{1-\nu(x,y)}}{\sqrt{\Ga(\nu(x))\Ga(\nu(y))}}M_{\nu(x,y)}\left(\sqrt{Q_{xy}(x-y)}\right),
  \end{equation}
where $\Sig_x$ is a positive definite matrix accounting for the local
anisotropy at location $x$, $\phi_{xy} =
|\Sig_x|^{1/4}|\Sig_y|^{1/4}\left|\frac{\Sig_x +
  \Sig_y}{2}\right|^{-1/2}$, $Q_{xy}(h) = h^t\left(\frac{\Sig_x +
  \Sig_y}{2}\right)^{-1}h$, $M_{\nu(x,y)}(h) =
h^{\nu(x,y)}K_{\nu(x,y)}(h)$, $K_{\nu(x,y)}$ is the modified Bessel
function of the second kind and $\nu(x,y) = \frac{\nu(x)+\nu(y)}{2}$.

A simulation method for this type of covariance has been recently
proposed in \citep{emerycontinuous}. The simulation is obtained as a
weighted sum of cosine waves, with random frequencies computed from an
instrumental stationary spectral density and random phases. A location
specific importance weight is applied so as to respect the local
spectral density of \eqref{eq:covmatns}.

\begin{figure}[h!]
  \centering
  \begin{subfigure}[t]{0.24\linewidth}
    \includegraphics[width=\textwidth]{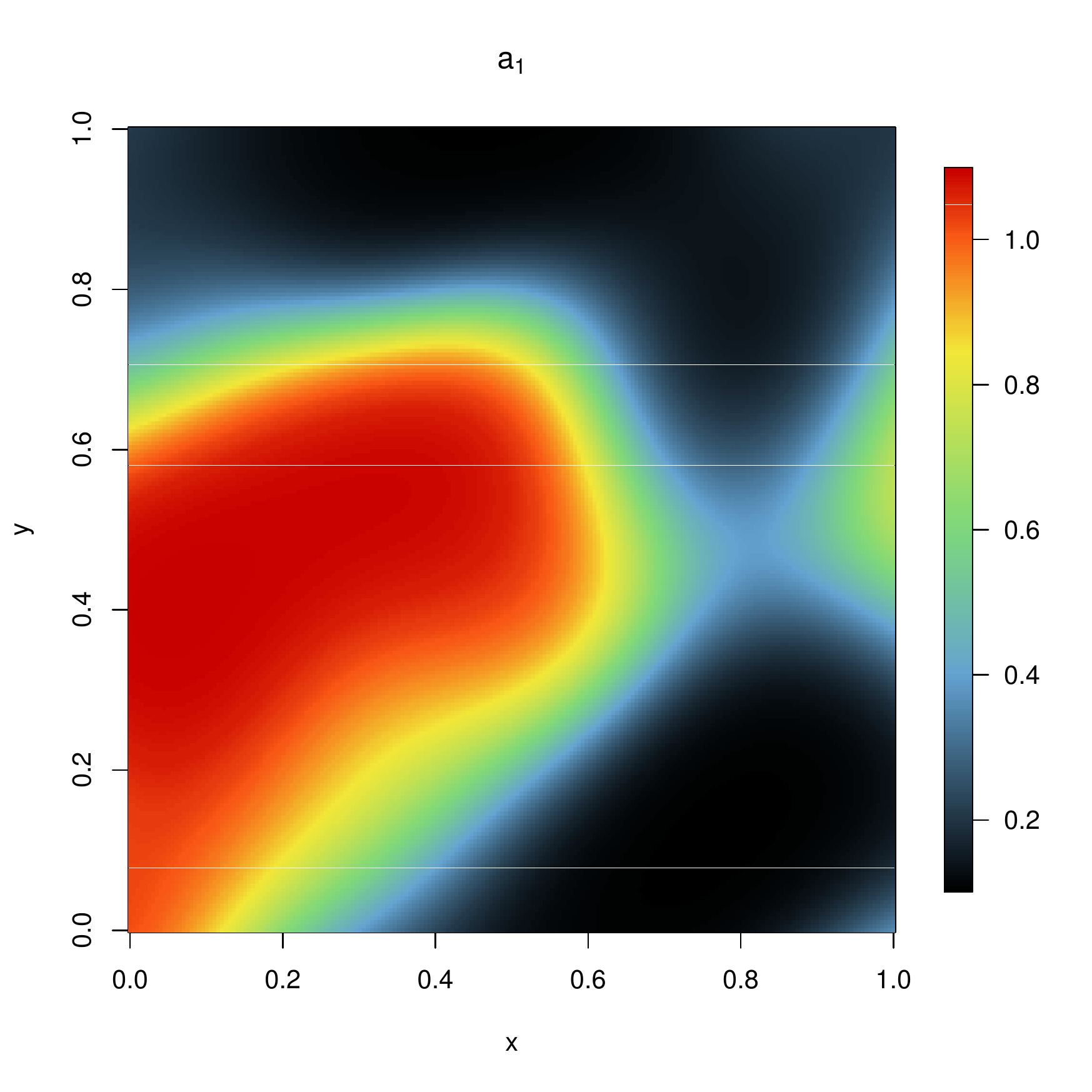}
    \caption{First scale parameter}
    \label{fig:simnsa1}
  \end{subfigure}
  \begin{subfigure}[t]{0.24\linewidth}
    \includegraphics[width=\textwidth]{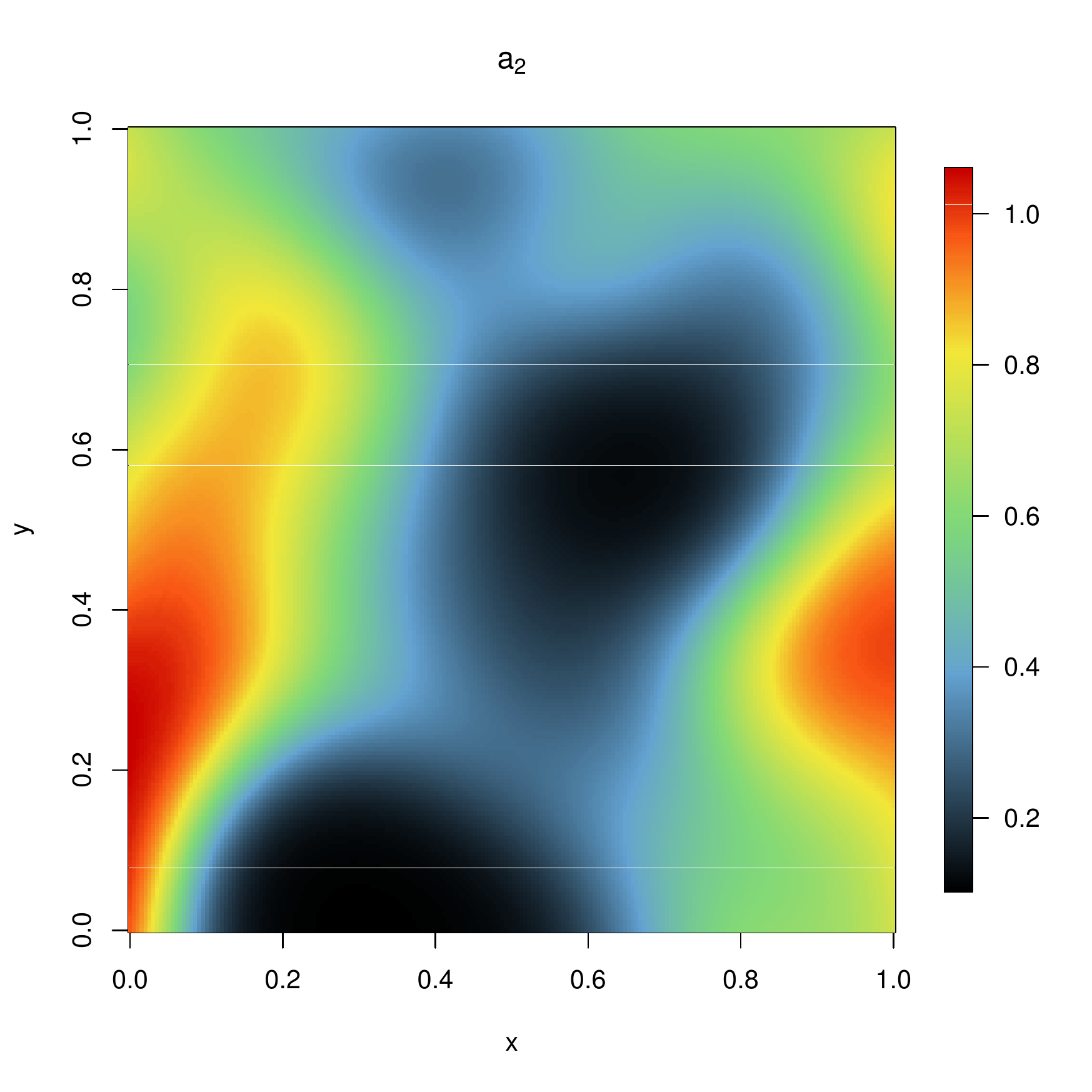}
    \caption{Second scale parameter}
    \label{fig:simnsa2}
  \end{subfigure}
  \begin{subfigure}[t]{0.24\linewidth}
    \includegraphics[width=\textwidth]{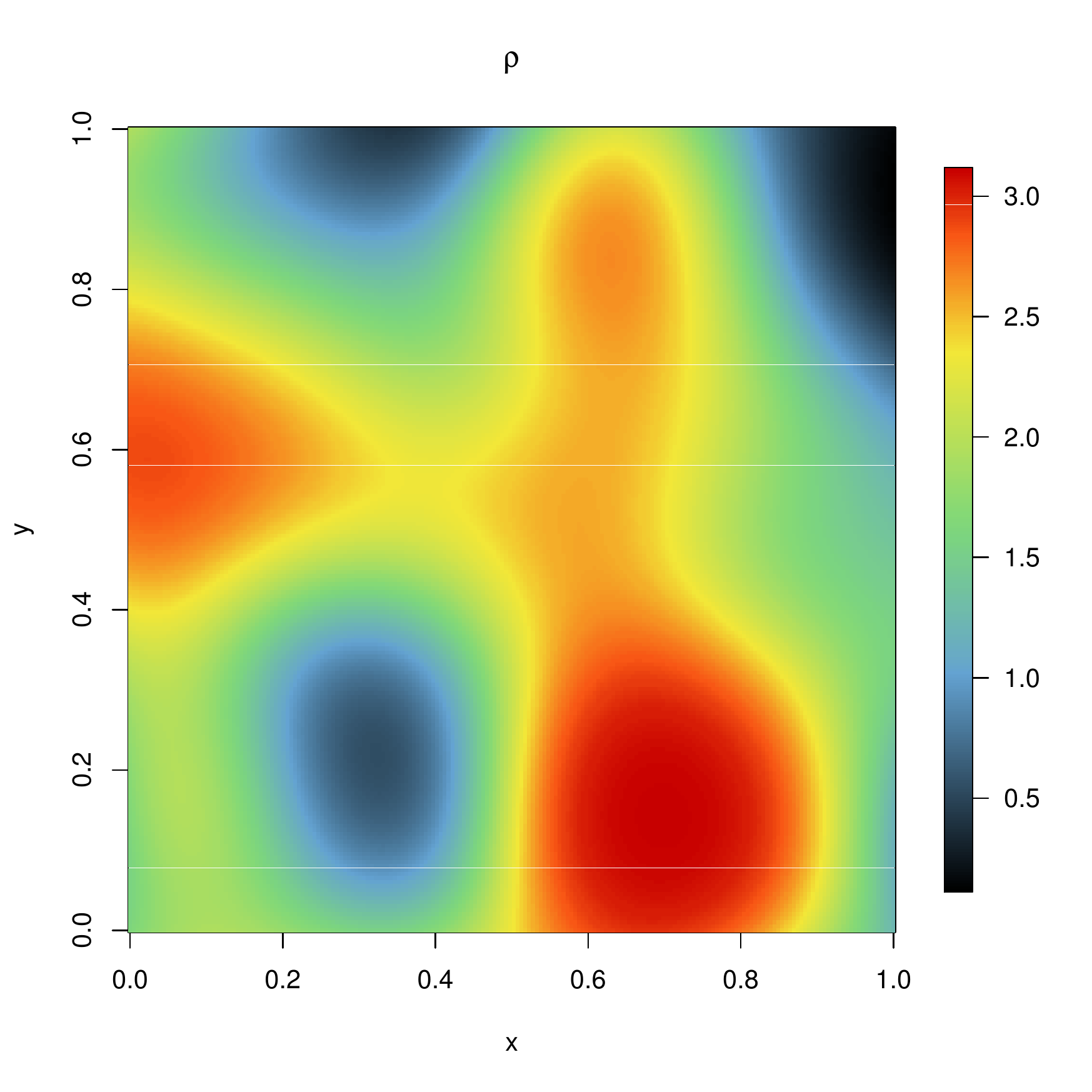}
    \caption{Anisotropy angle}
    \label{fig:simnstet}
  \end{subfigure}
  \begin{subfigure}[t]{0.24\linewidth}
    \includegraphics[width=\textwidth]{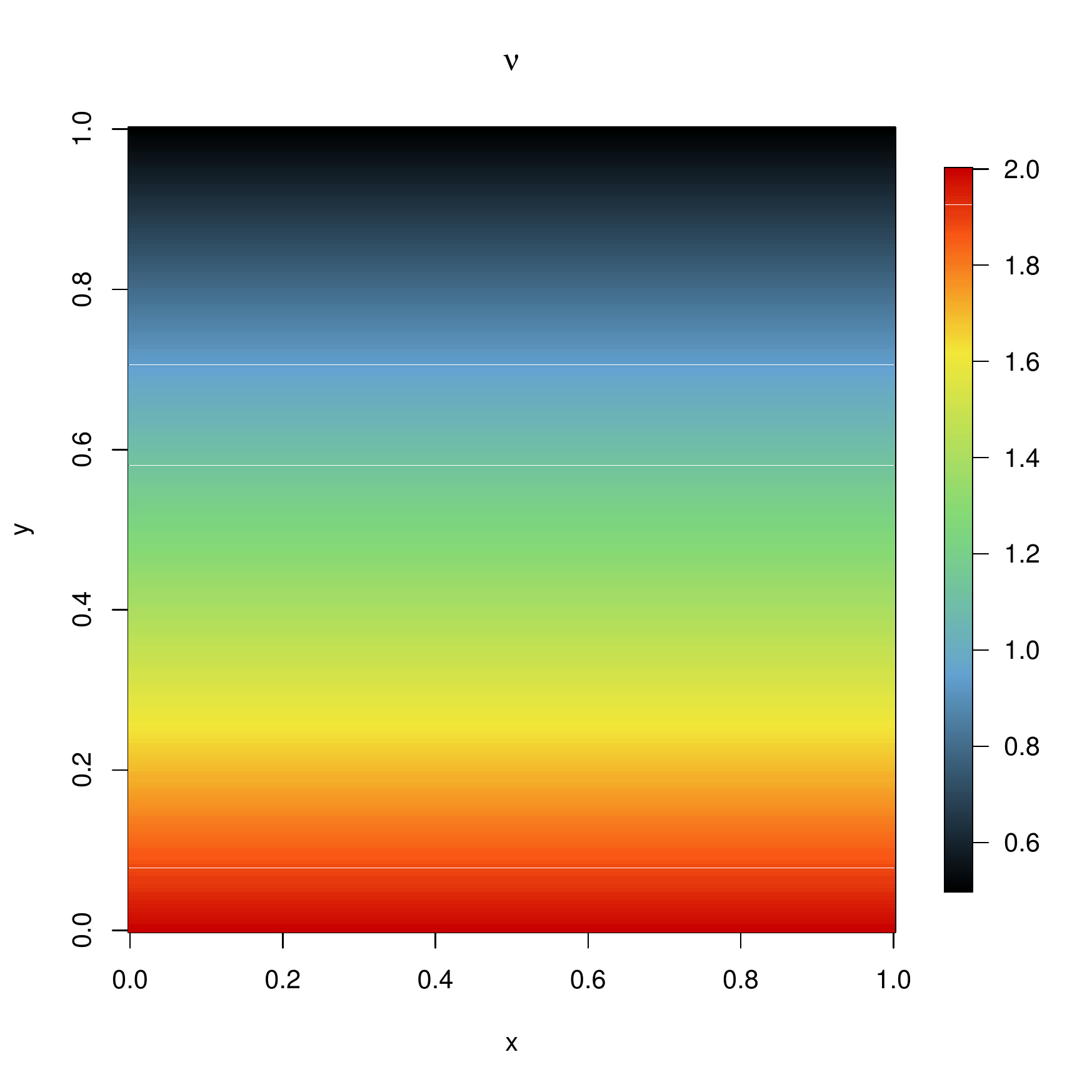}
    \caption{Smoothness parameter}
    \label{fig:simnsnu}
  \end{subfigure}
  \caption{Covariance parameters}
\label{fig:param_ns}
\end{figure}

To build our example, we first simulate a map of varying parameters
over the unit square. The two scale parameters and the anisotropy
angle are simulated as uniform transforms of a Gaussian random field
with a gaussian covariance and a large scale parameter. The scale
parameters are comprised between 0.1 and 1.1, while the anisotropy
angle is allowed to vary between 0 and $\pi$. The smoothness parameter
varies linearly between 0.5 and 2 along the vertical direction. The
varying parameters are plotted in figure \ref{fig:param_ns}.

\begin{figure}[h!]
  \centering
  \begin{subfigure}[t]{0.45\linewidth}
    \includegraphics[width=\textwidth]{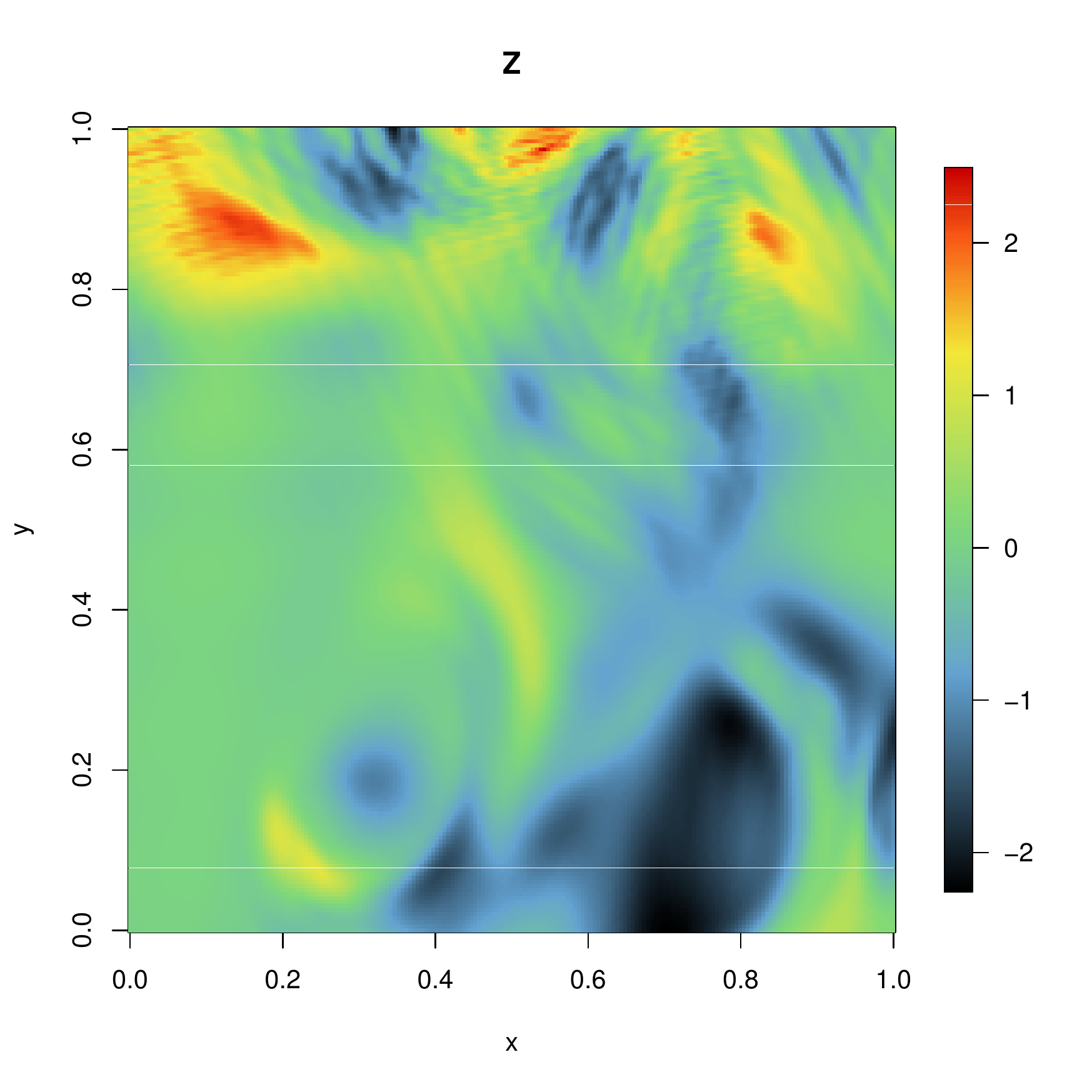}
    \caption{Reference realization}
    \label{fig:simref_ns}
  \end{subfigure}
  \begin{subfigure}[t]{0.45\linewidth}
    \includegraphics[width=\textwidth]{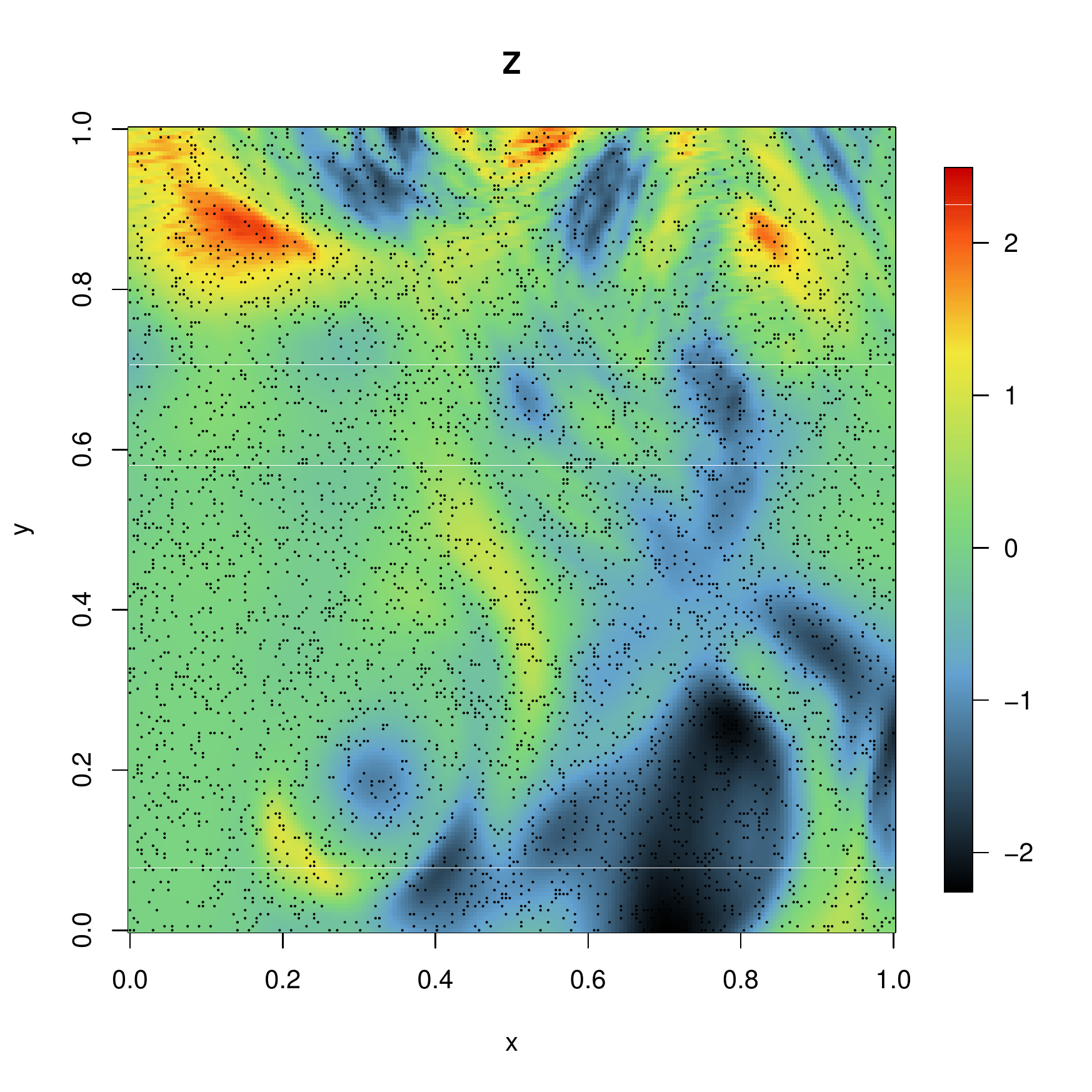}
    \caption{Sampling at 5000 locations}
    \label{fig:simsam_ns}
  \end{subfigure}
\caption{Reference simulation and sampling}
\label{fig:sim_ns}
\end{figure}

10000 cosine waves are used to generate the realisation plotted in
figure \ref{fig:sim_ns}. We can clearly see the effect of the varying
anisotropy and smoothness parameters. Besides, some artifacts appear
at the top of the domain, where the smoothness parameter is the
lowest. They are due to a too small number of cosines in the
simulation, which implies a incomplete sampling of the highest
frequencies. This is not an important problem in this example however
as these artifacts will be treated as patterns in the original
image. This realisation is sampled at 5000 locations uniformly
distributed over the domain, as in the previous example.

The inference is conducted following the same steps as in the previous
example. Moreover, the same set of 2658 basis functions is used to
constitute the dictionary of candidates.

\begin{figure}[h!]
  \centering
  \begin{subfigure}[t]{0.45\linewidth}
    \includegraphics[width=\textwidth]{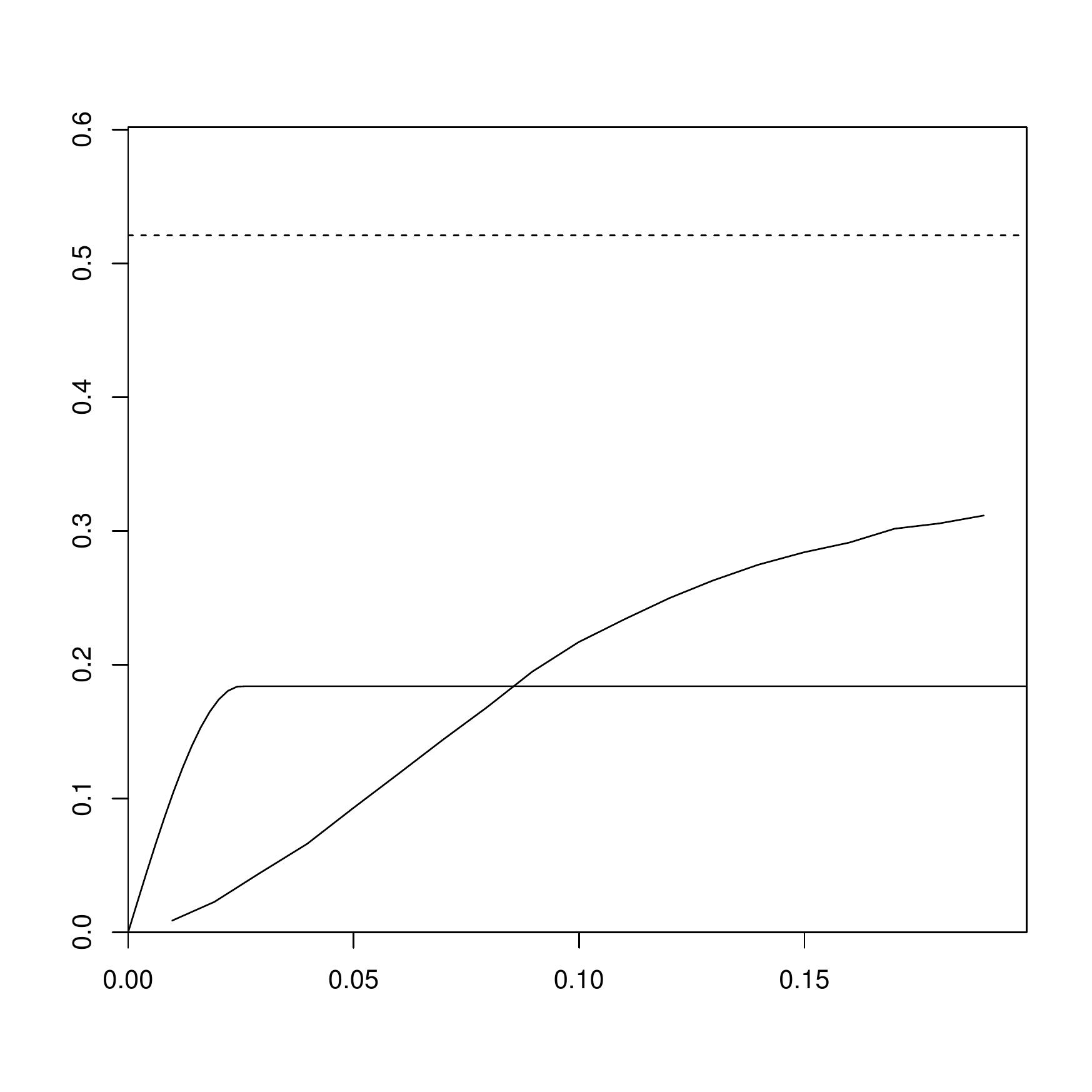}
    \caption{First fitting of the small scale structure}
    \label{fig:varions1}
  \end{subfigure}
  \begin{subfigure}[t]{0.45\linewidth}
    \includegraphics[width=\textwidth]{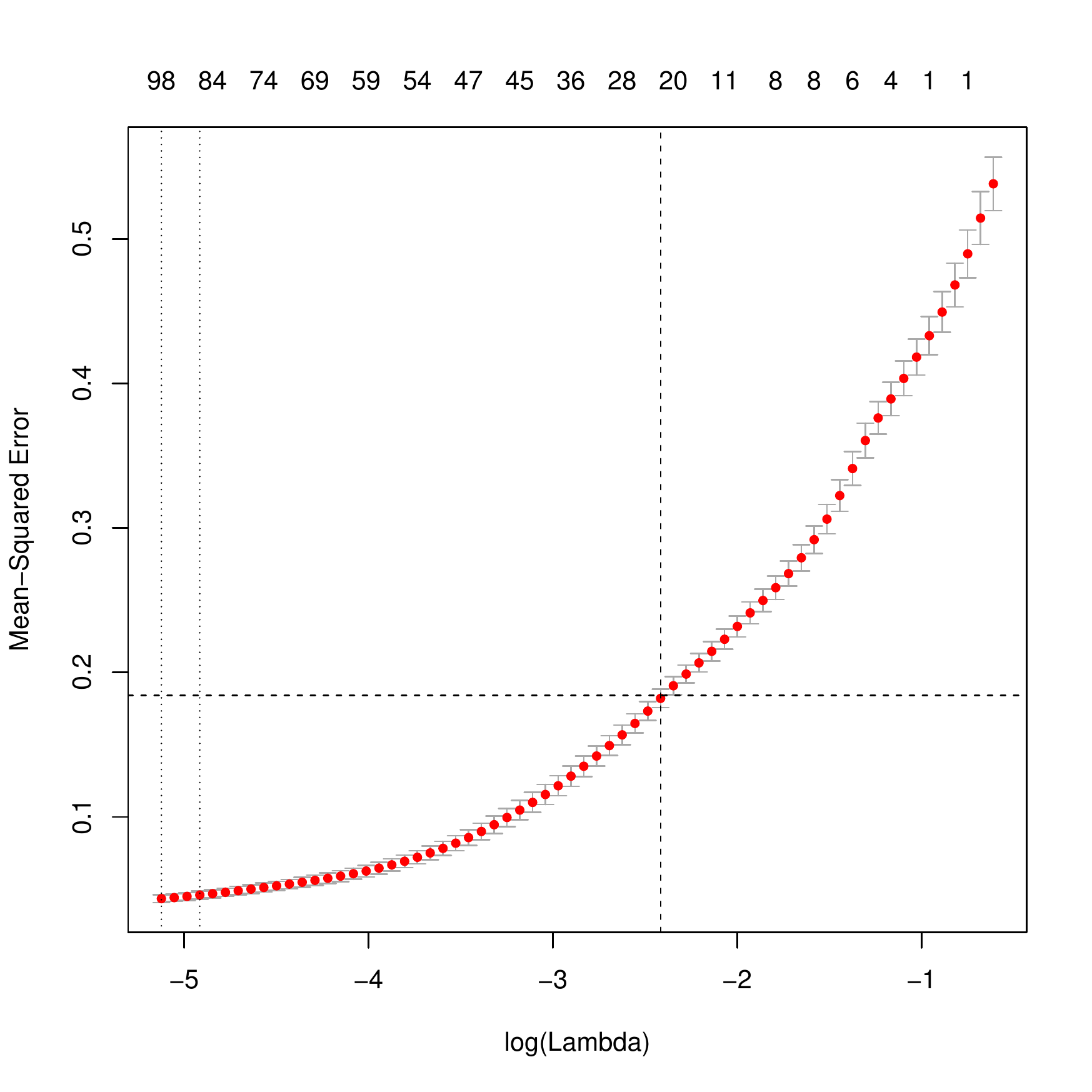}
    \caption{First selection of basis functions}
    \label{fig:glmns1}
  \end{subfigure}
  \caption{First fitting of the model}
\end{figure}

We can see in figure \ref{fig:varions1} that the variogram fit is not
good, which is due to the low value for the range of the taper
compared to what should be the range. In particular, a too large
weight is associated to the small scale component, which results in
the selection of only 23 basis functions.

\begin{figure}[h!]
  \centering
  \begin{subfigure}[t]{0.45\linewidth}
    \includegraphics[width=\textwidth]{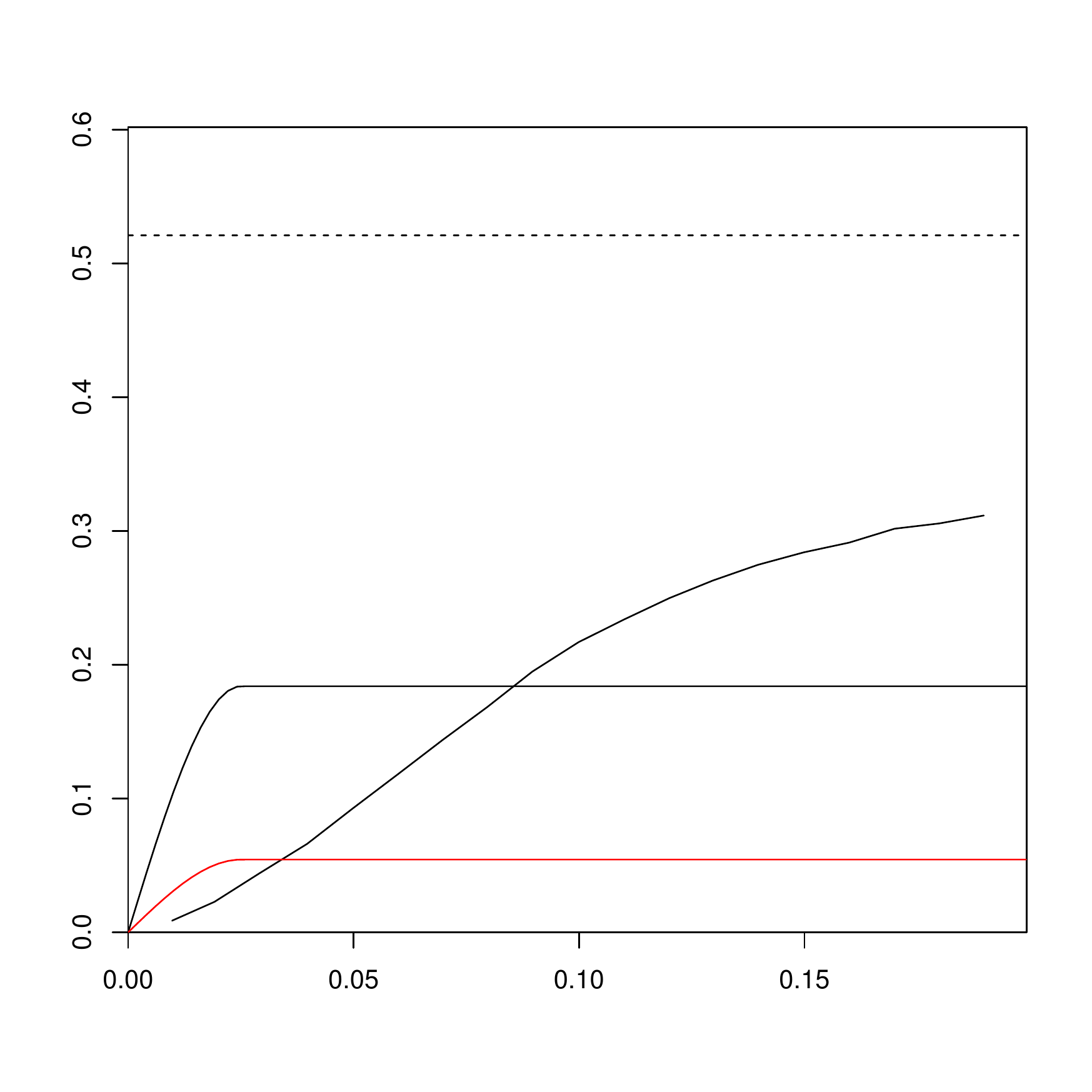}
    \caption{Updated variogram}
    \label{fig:varions2}
  \end{subfigure}
  \begin{subfigure}[t]{0.45\linewidth}
    \includegraphics[width=\textwidth]{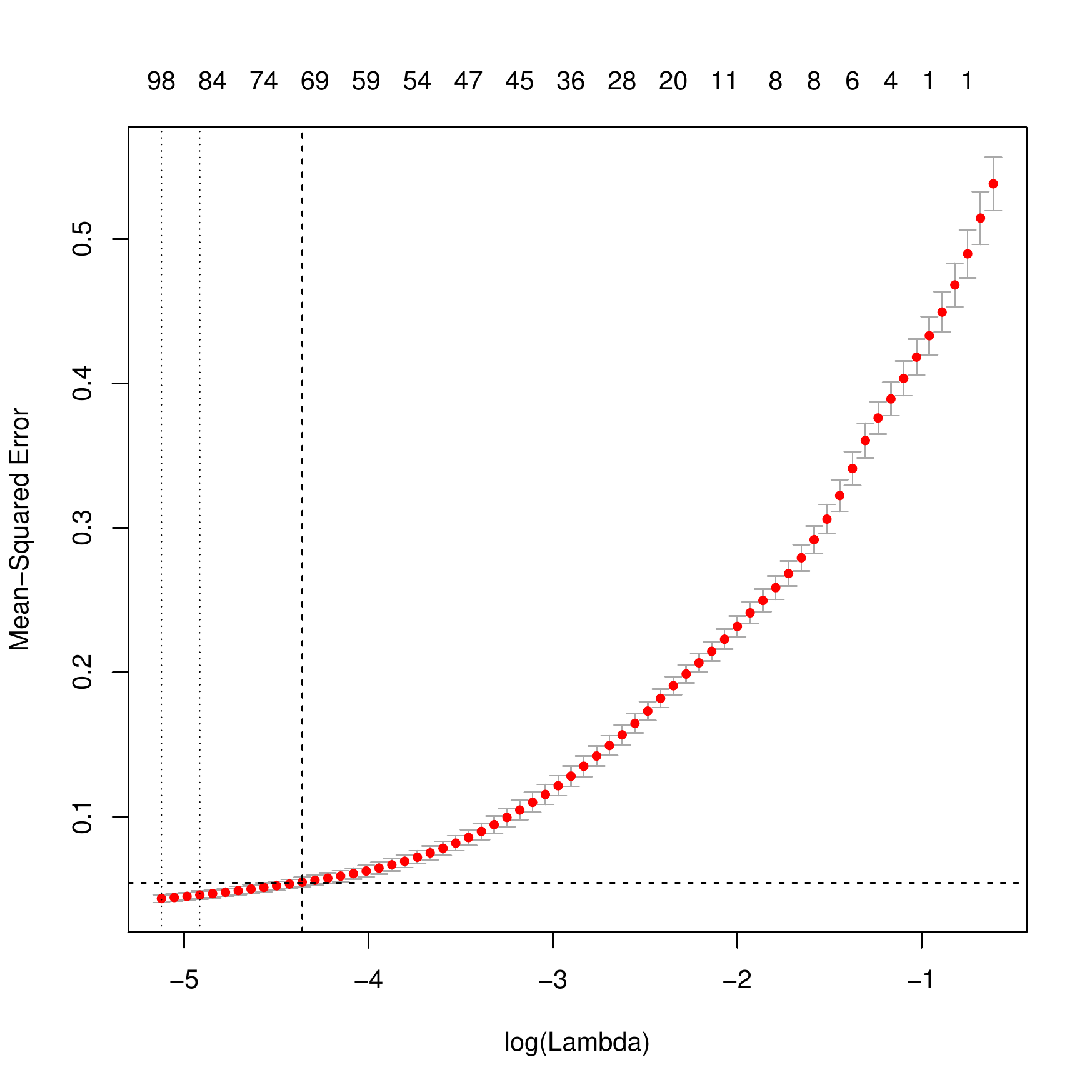}
    \caption{Updated selection of basis functions}
    \label{fig:glmns2}
  \end{subfigure}
  \caption{Model update}
\end{figure}

After the EM algorithm, the updated value of $\sig^2$ provides a
better fit of the variogram and 69 basis functions are now selected.

\begin{figure}[h!]
  \centering
  \begin{subfigure}[t]{0.45\linewidth}
    \includegraphics[width=\textwidth]{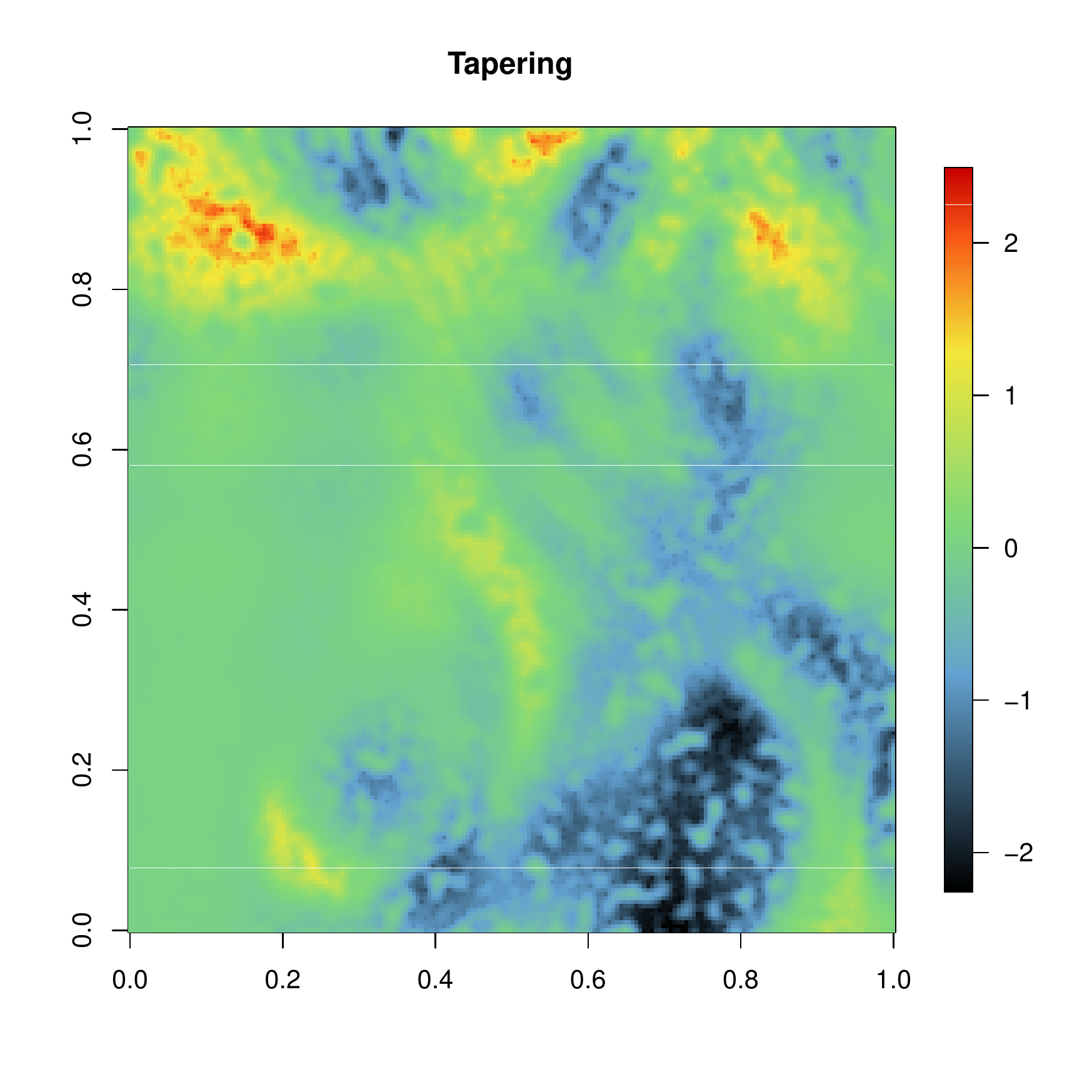}
    \caption{Covariance tapering}
    \label{fig:krigns_tap}
  \end{subfigure}
  \begin{subfigure}[t]{0.45\linewidth}
    \includegraphics[width=\textwidth]{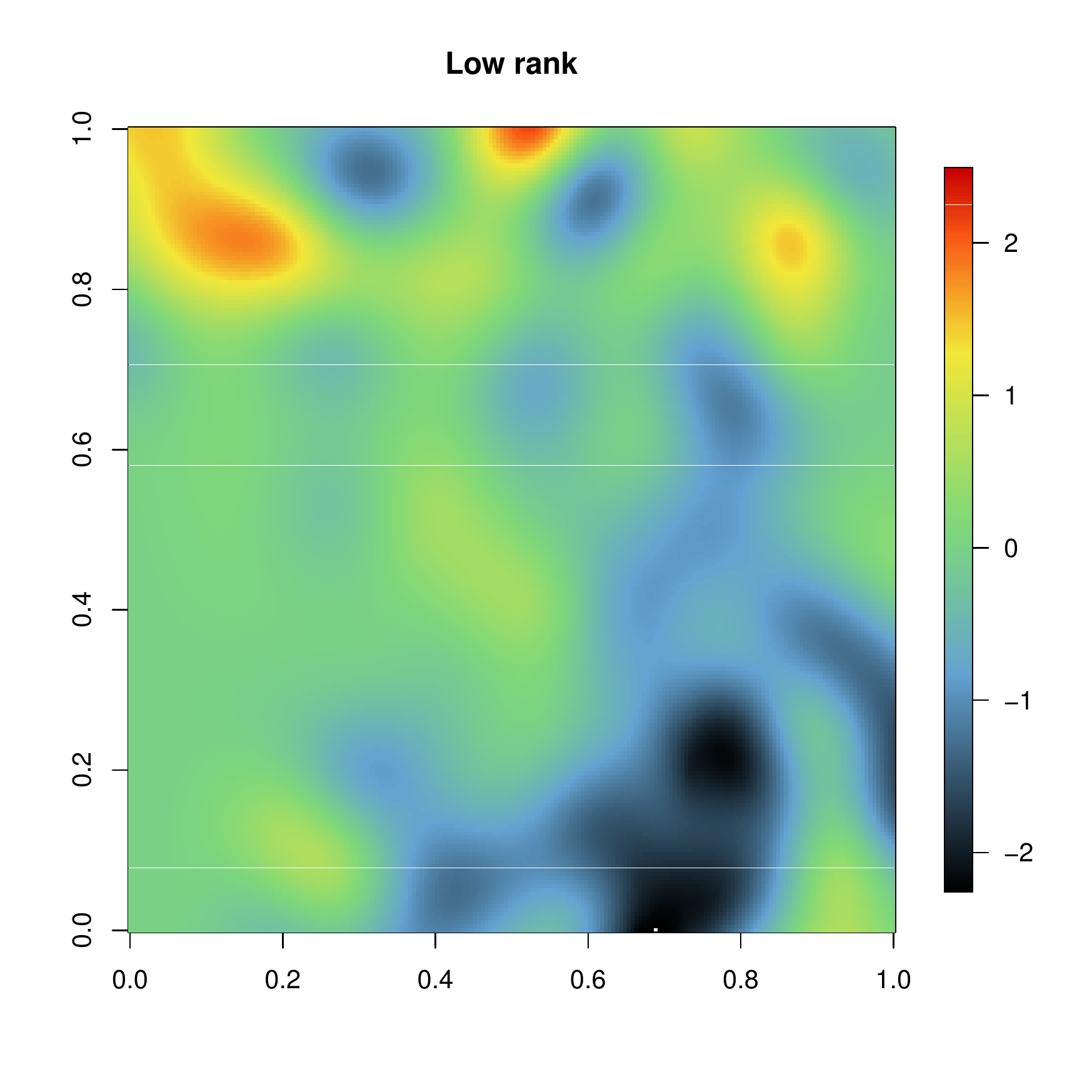}
    \caption{Low rank model}
    \label{fig:krigns_lr}
  \end{subfigure}
  \begin{subfigure}[t]{0.45\linewidth}
    \includegraphics[width=\textwidth]{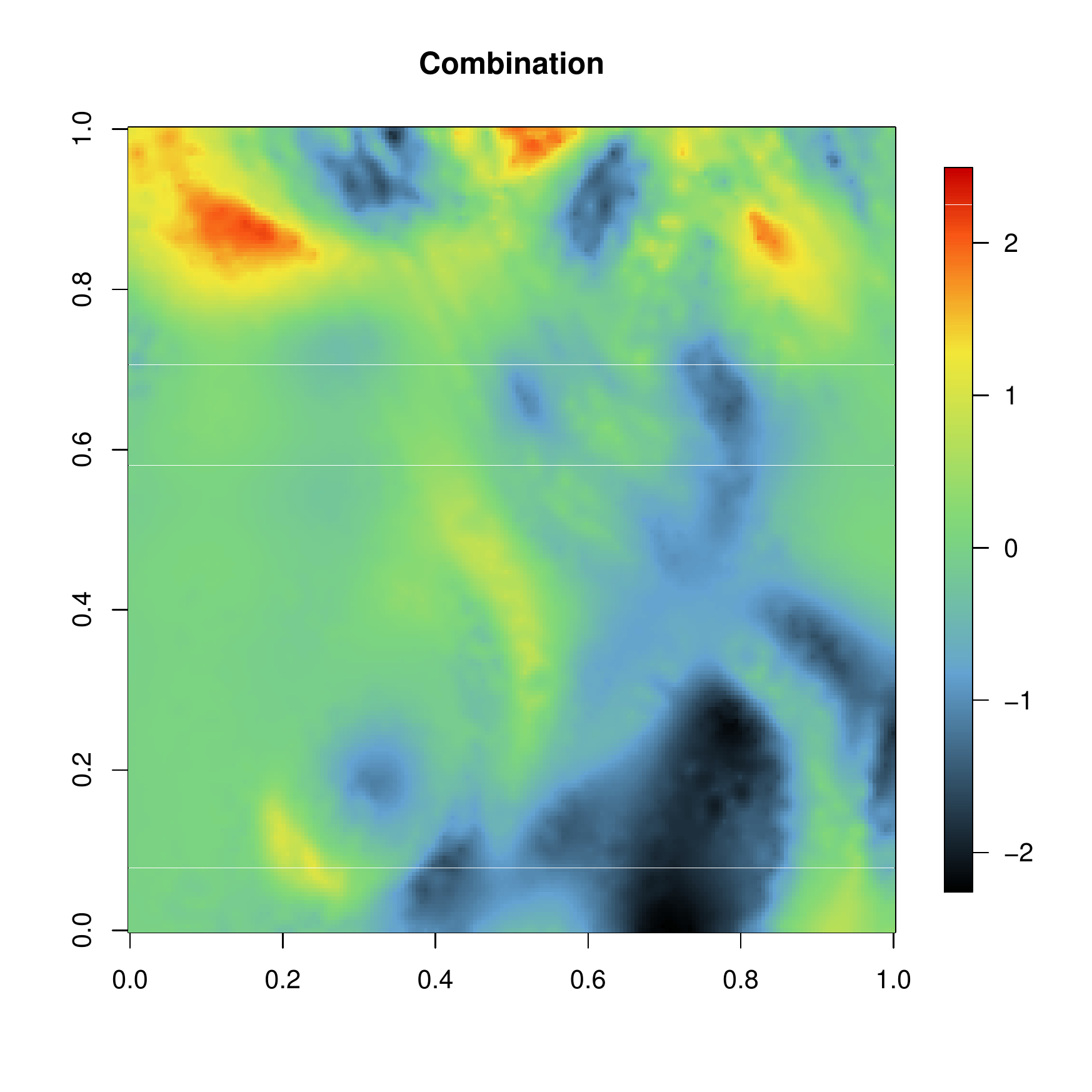}
    \caption{Combination}
    \label{fig:krigns_combi}
  \end{subfigure}
  \begin{subfigure}[t]{0.45\linewidth}
    \includegraphics[width=\textwidth]{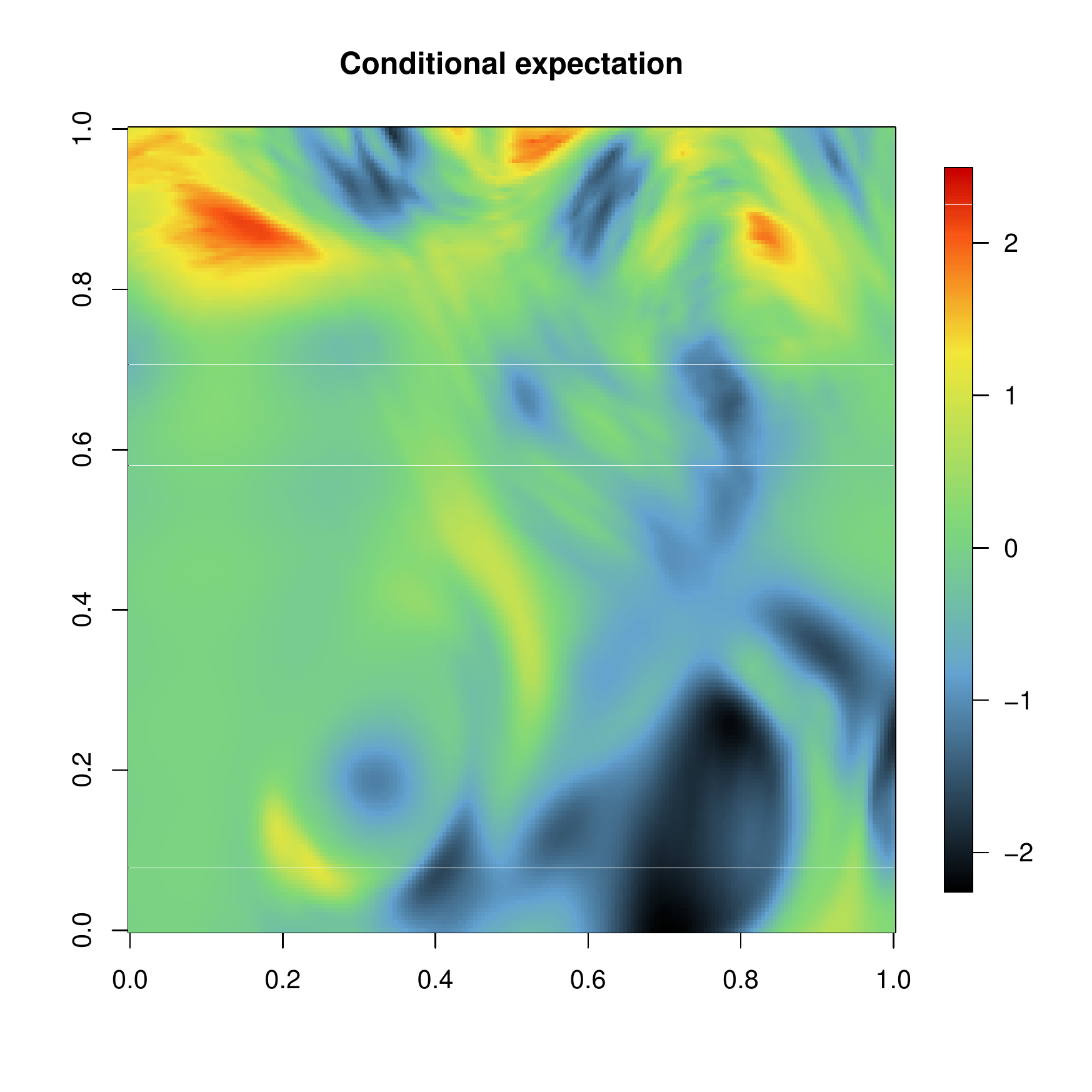}
    \caption{Conditional expectation}
    \label{fig:krigns_ce}
  \end{subfigure}
\caption{Kriging maps}
\label{fig:krigns_maps}
\end{figure}

The difference in the kriging maps shown in figure
\ref{fig:krigns_maps} are clearly visible. The small value chosen for
the tapering range causes the prediction to be equal to the mean when
the distance from observation locations increase. The low rank model
provides on the contrary an oversmooth map. The combination allows
capturing both scales of variation and offers a result closer to the
conditional expectation. In terms of mean square prediction error
(MSPE), an important improvement is obtained with the combination as
can be read in table \ref{tab:mspens}.

\begin{figure}[h!]
  \centering
  \begin{subfigure}[t]{0.3\linewidth}
    \includegraphics[width=\textwidth]{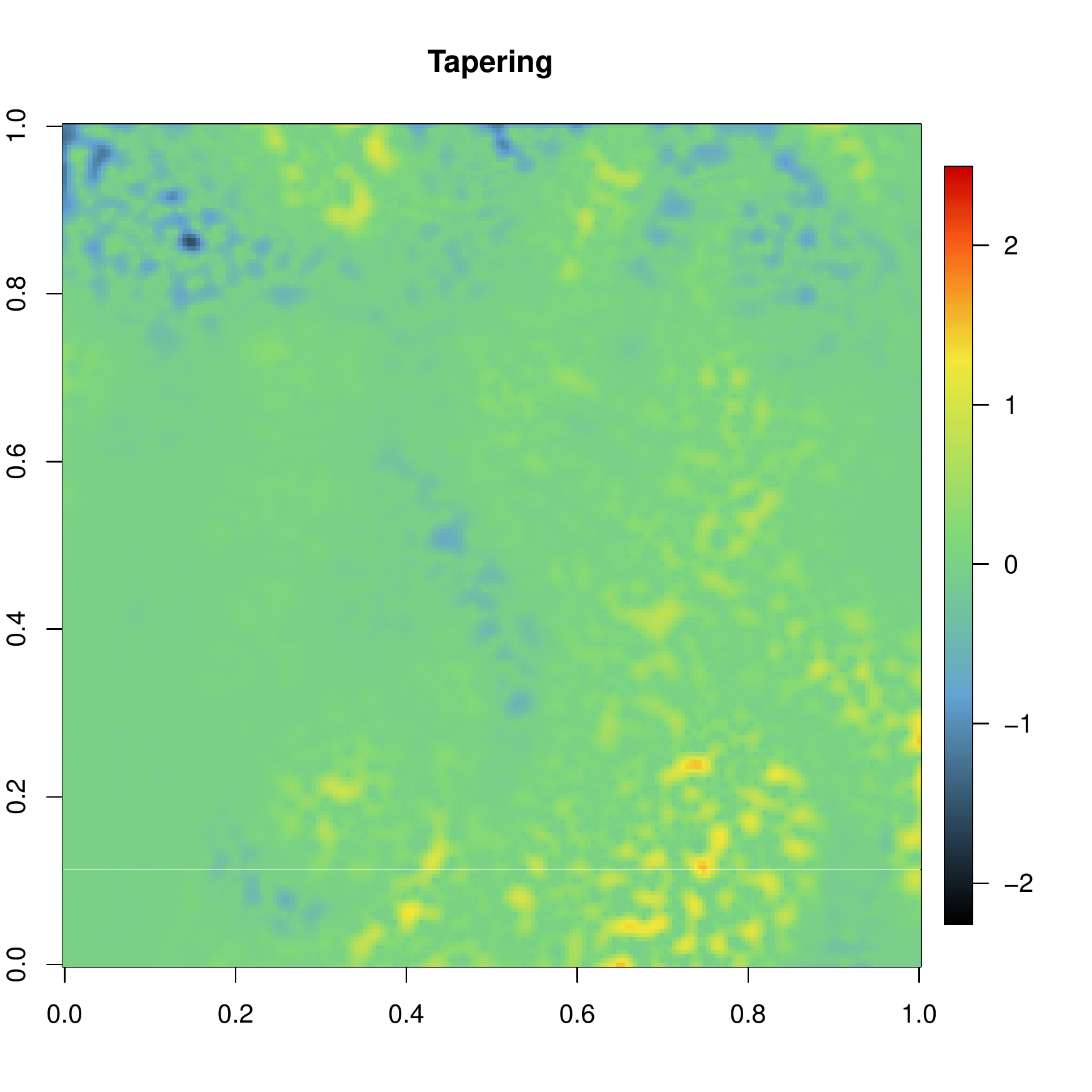}
    \caption{Covariance tapering}
    \label{fig:krignsdiff_tap}
  \end{subfigure}
  \begin{subfigure}[t]{0.3\linewidth}
    \includegraphics[width=\textwidth]{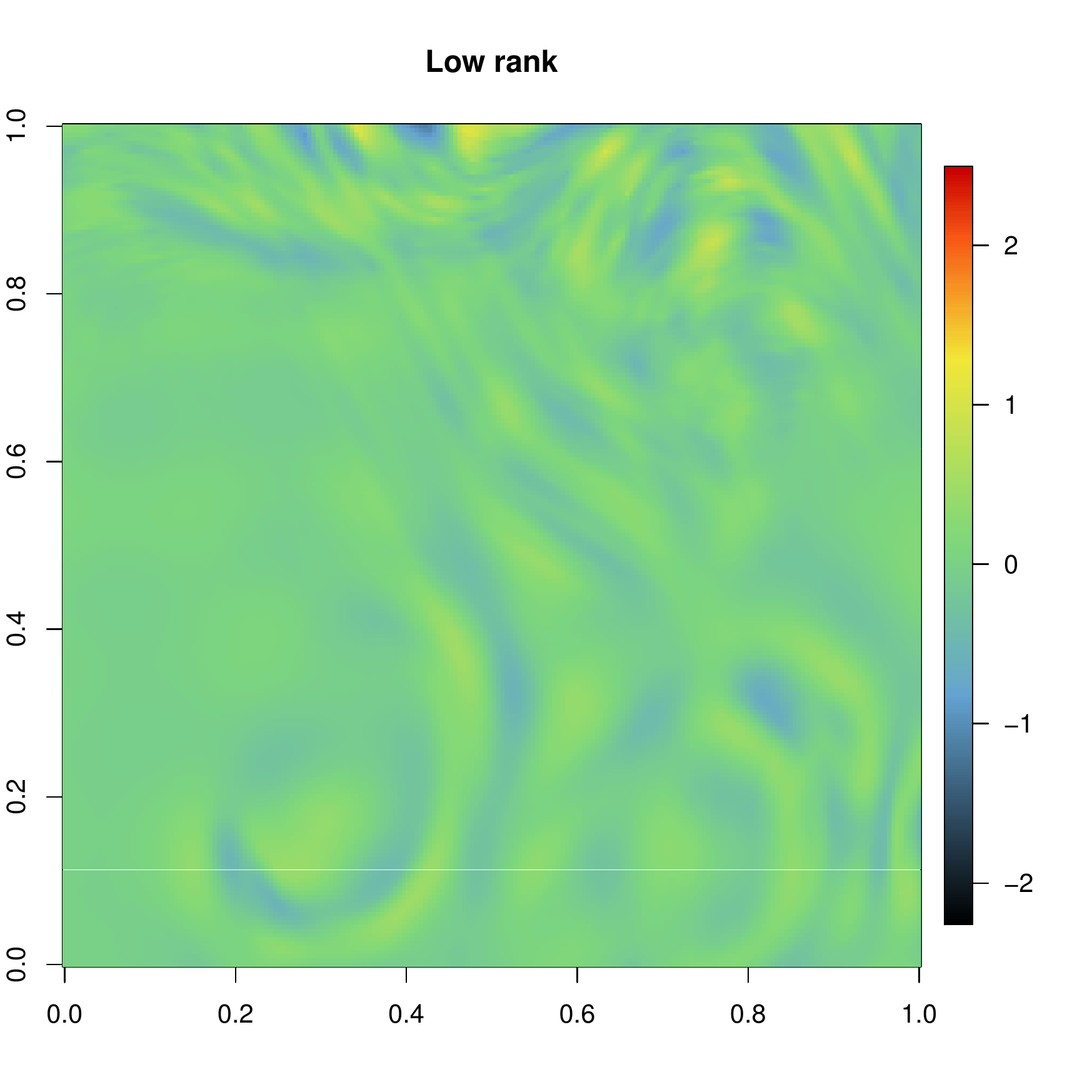}
    \caption{Low rank model}
    \label{fig:krignsdiff_lr}
  \end{subfigure}
  \begin{subfigure}[t]{0.3\linewidth}
    \includegraphics[width=\textwidth]{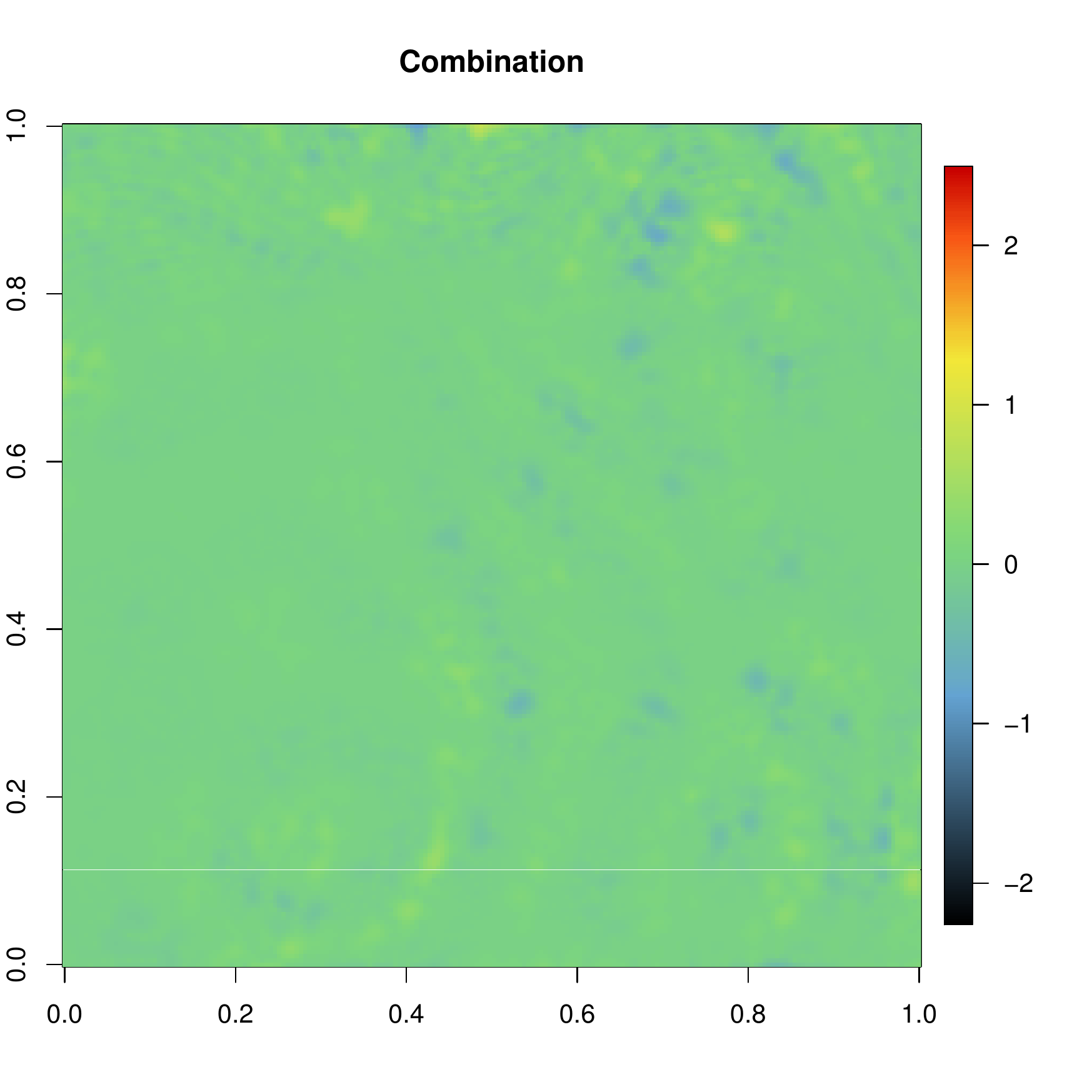}
    \caption{Combination}
    \label{fig:krignsdiff_combi}
  \end{subfigure}
\caption{Difference maps}
\label{fig:krigns_diff}
\end{figure}

\begin{table}
  \caption{Summary}
  \begin{tabular}{lccc}
    \hline\\
    model & small scale structure & number of basis functions & MSPE \\
    \hline\\

    Covariance tapering & range 0.025 & 0 & 0.048 \\

    Low rank  & nugget 0.03 & 560 & 0.036 \\

    Combination & range 0.025 & 287 & 0.008\\

    Conditional expectation & & & 0.002\\
    \hline
  \end{tabular}
  \label{tab:mspens}
\end{table}

\subsection{Total column ozone dataset}
We now apply our approach to total column ozone (TCO) data acquired
from an orbiting satellite mounted with a passive sensor registering
backscattered solar ultraviolet radiation. The dataset we consider
here has been previously analyzed in \citet{fixedrank},
\citet{bolin2011spatial} and \citet{eidsvik2014estimation} among
others. It consists of $n$ = 173 405 measurements.

We compare covariance tapering, fixed rank kriging (FRK) and the
proposed combination approach.

The predicition is performed on a 180 $\f$ 288 grid. The latitude
ranges from -89.5 to 89.5 in 1$^{\circ}$ steps, while longitude ranges from
-179.375 to 179.375 in 1.25$^{\circ}$ steps. This represents 51 840 prediction
sites in total.

Covariance tapering is implemented using a Wendland2 \citep{compcov}
covariance with a taper range set at 275 km, which represents a
2.5$^{\circ}$ distance. Indeed, the data exhibit a rather smooth
behavior that can be well represented by a regular covarariance
function. The parameter estimation is performed through variogram
fitting. FRK is performed with 4 resolution bi-square basis functions,
representing a total of 3608 basis functions. The combination is built
using the same fine scale structure as in covariance tapering and by
selecting the basis functions among those used in FRK. 331 basis
functions are selected. The EM step is run only once since the first
fit of the variogram is deemed satisfying.

\begin{figure}[h!]
  \centering
  \begin{subfigure}[t]{0.49\linewidth}
    \includegraphics[width=\textwidth]{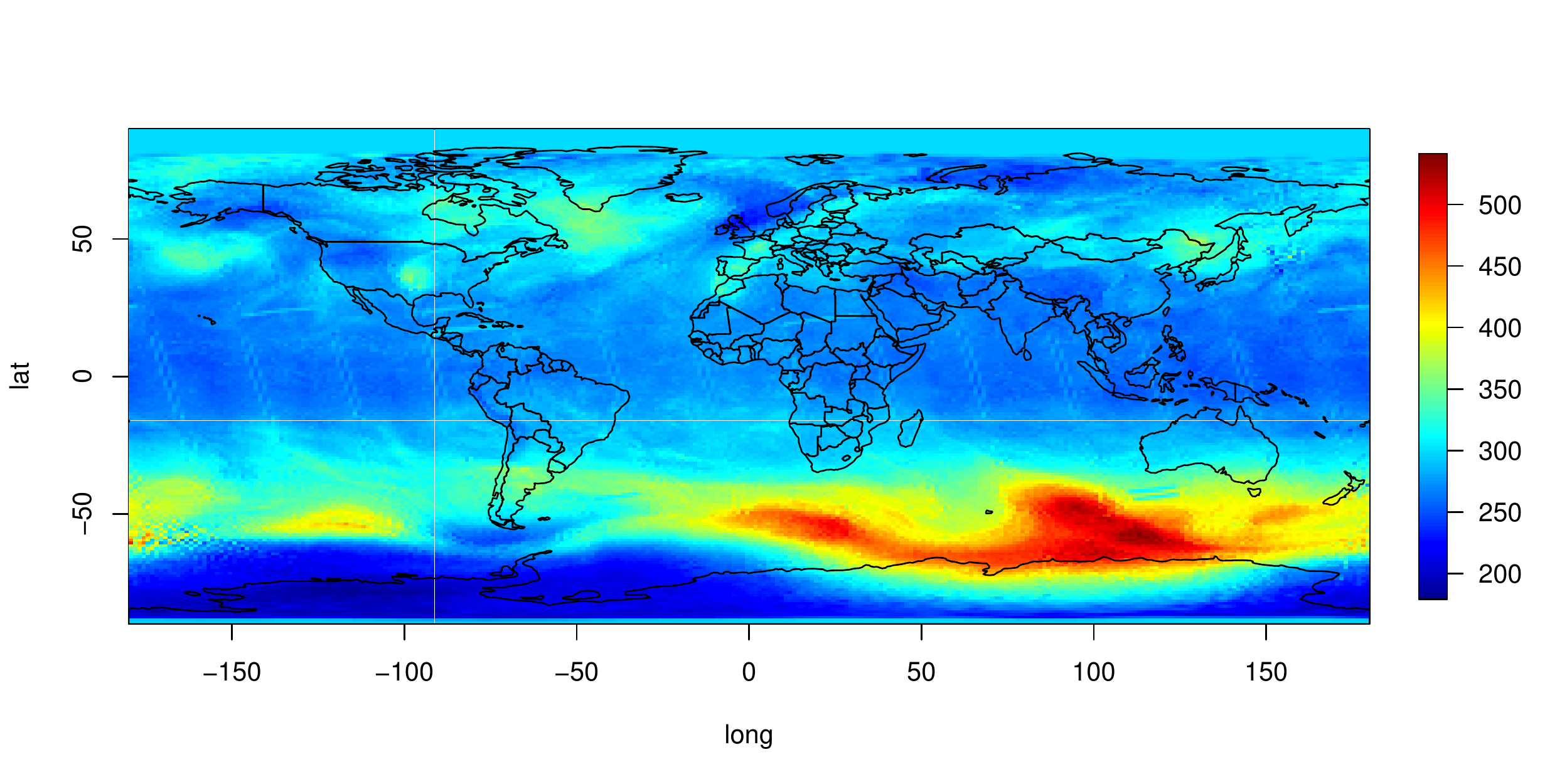}
    \caption{Covariance tapering prediction}
    \label{fig:tomstap}
  \end{subfigure}
  \begin{subfigure}[t]{0.49\linewidth}
    \includegraphics[width=\textwidth]{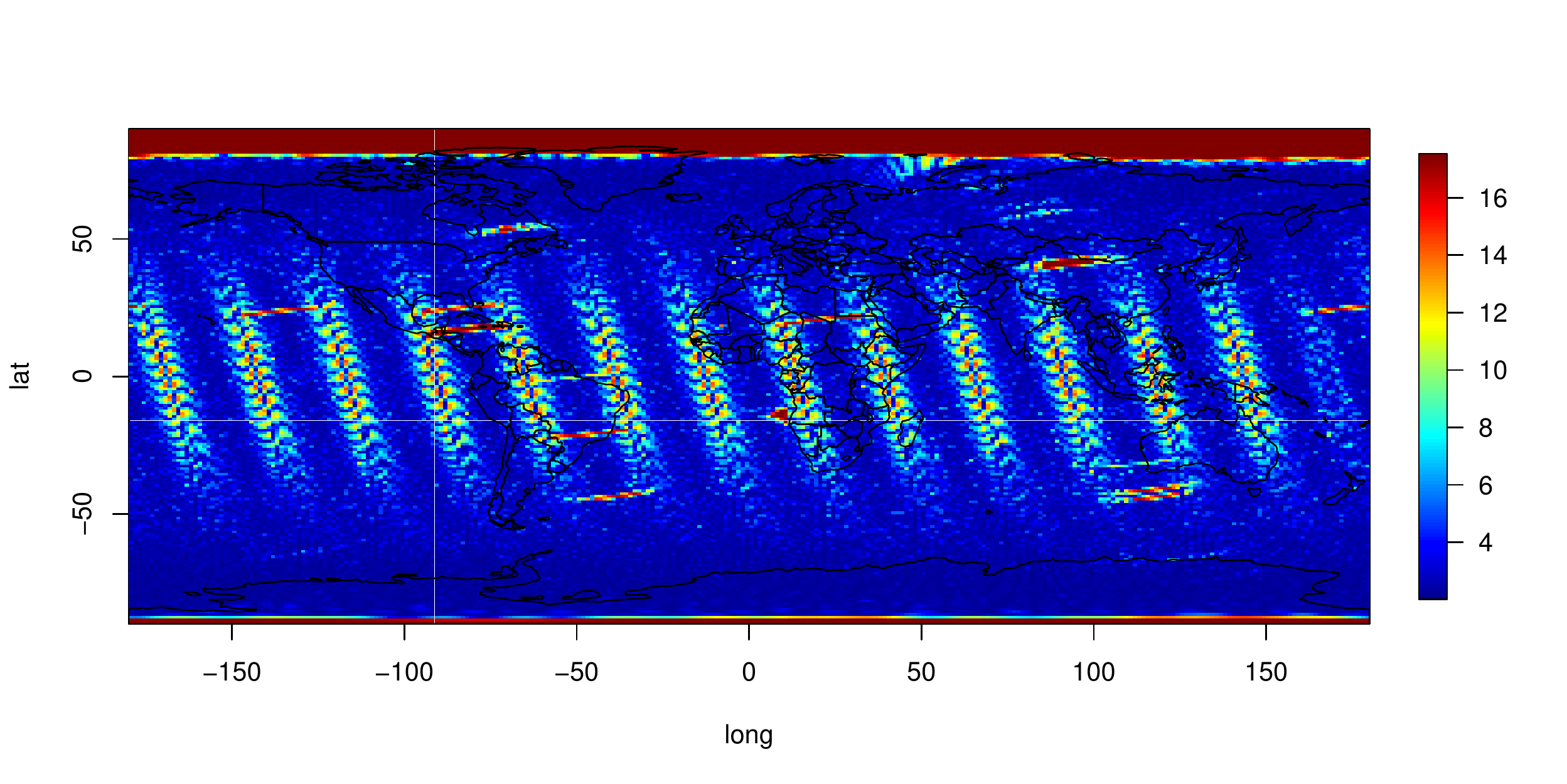}
    \caption{Covariance tapering standard deviation}
    \label{fig:tomstaps}
  \end{subfigure}
  \begin{subfigure}[t]{0.49\linewidth}
    \includegraphics[width=\textwidth]{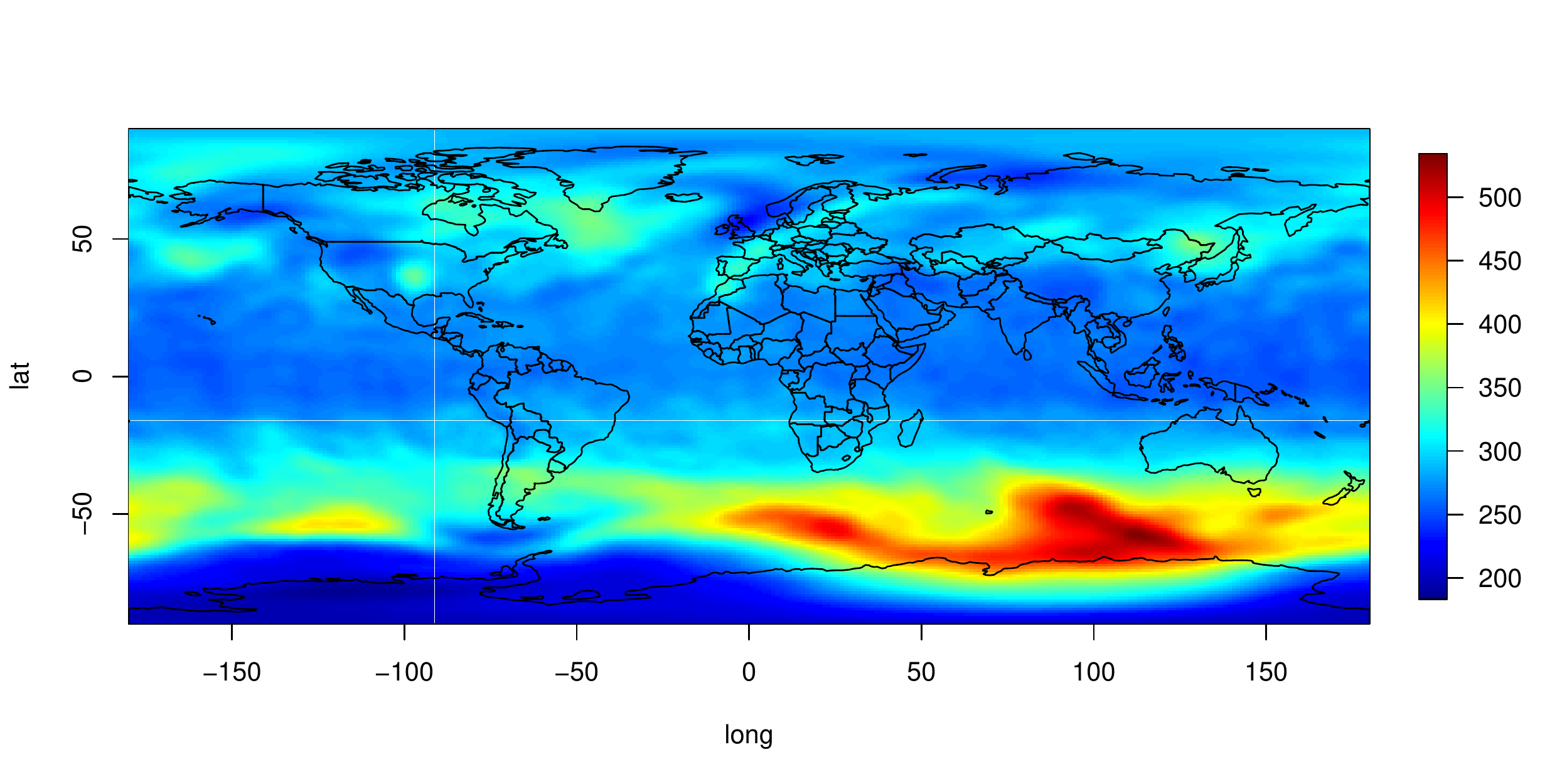}
    \caption{FRK prediction}
    \label{fig:tomsfrk}
  \end{subfigure}
  \begin{subfigure}[t]{0.49\linewidth}
    \includegraphics[width=\textwidth]{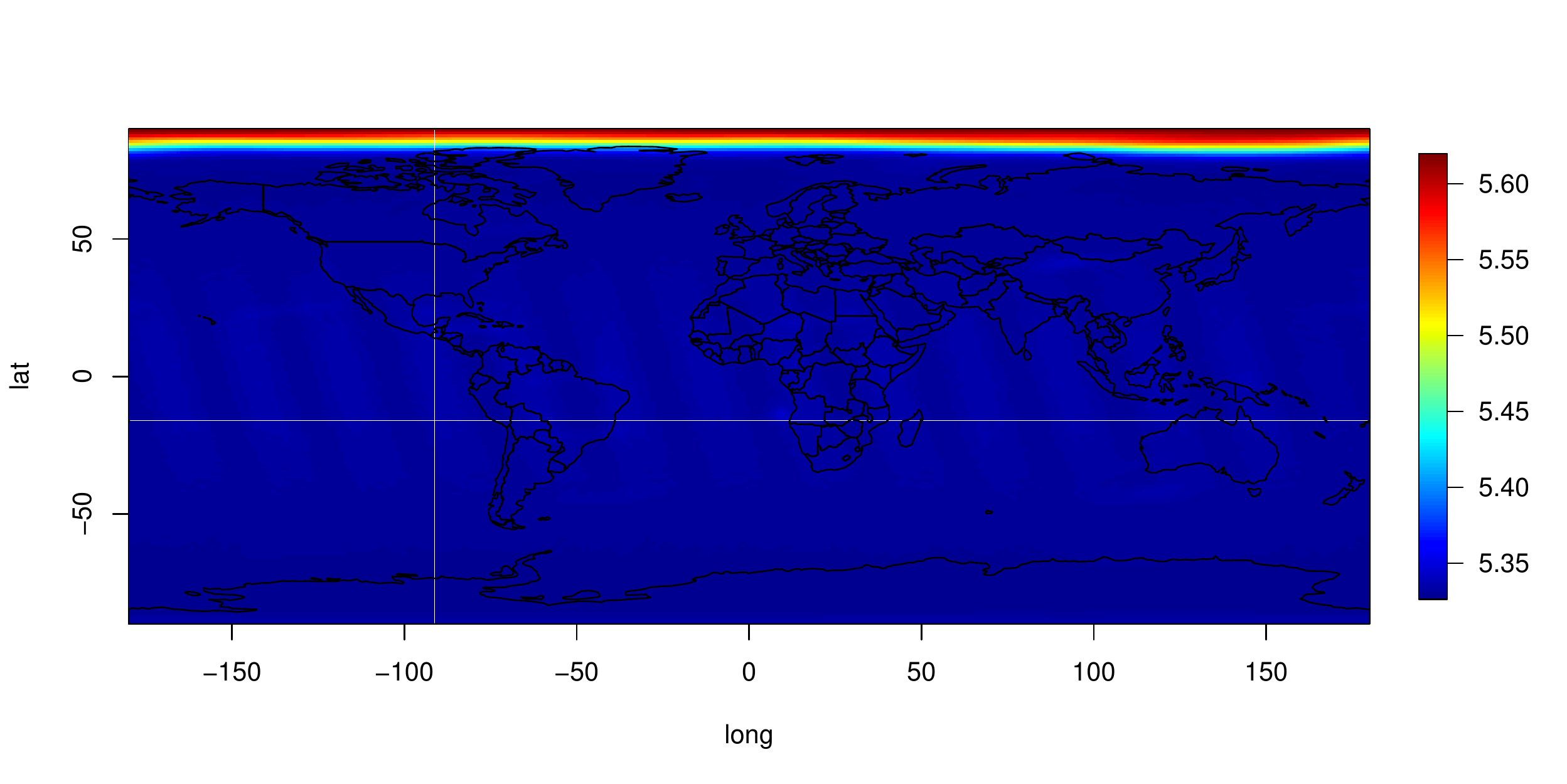}
    \caption{FRK standard deviation}
    \label{fig:tomsfrks}
  \end{subfigure}
    \begin{subfigure}[t]{0.49\linewidth}
    \includegraphics[width=\textwidth]{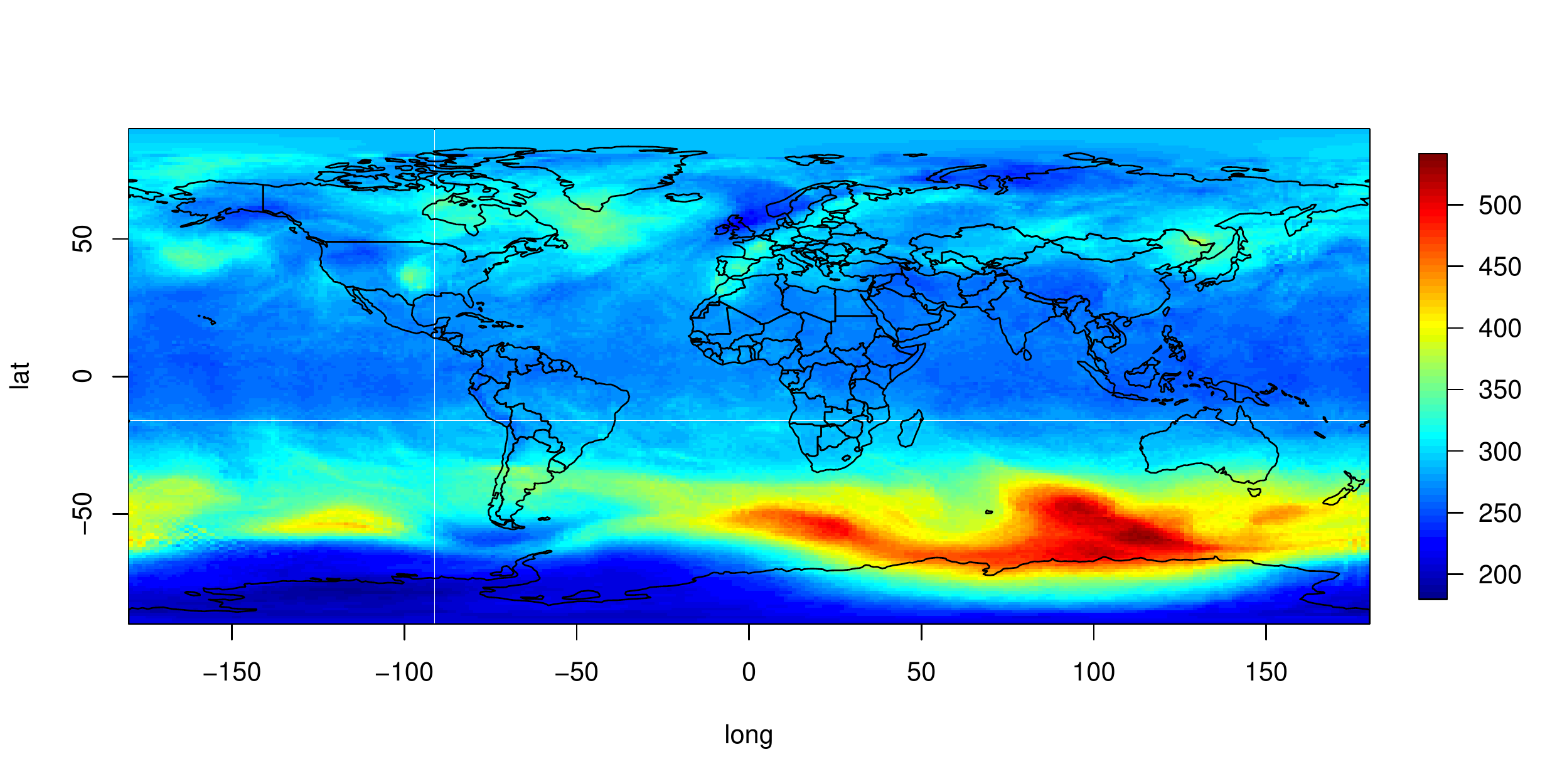}
    \caption{Combination prediction}
    \label{fig:tomscmb}
  \end{subfigure}
  \begin{subfigure}[t]{0.49\linewidth}
    \includegraphics[width=\textwidth]{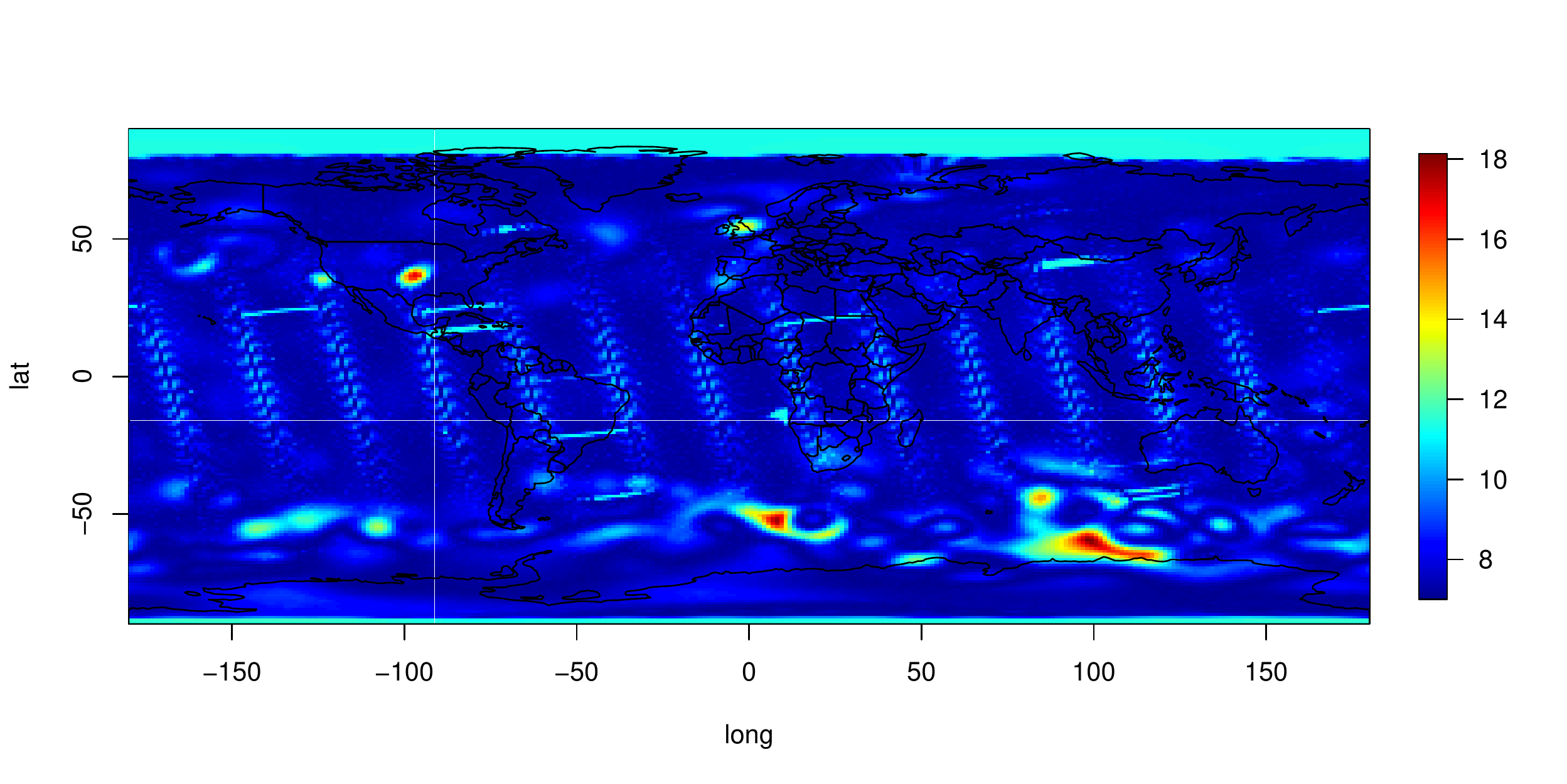}
    \caption{Combination standard deviation}
    \label{fig:tomscmbs}
  \end{subfigure}

\caption{Total Column Ozone data: Prediction and prediction standard
    deviation maps.}
\label{fig:toms}
\end{figure}

Figure \ref{fig:toms} shows the prediction maps of TCO and the
prediction standard deviations obtained by each three methods.  We can
observe that FRK predictions (fig. \ref{fig:tomsfrk}) are much
smoother than combination (fig. \ref{fig:tomscmb}) and tapering
(fig. \ref{fig:tomstap}) predictions, where more details are
visible. The latter presents however artifacts where the satellite
coverage is less dense, especially visible along lines going
south-southwest all around the globe. This effect is due to the use of
a too small tapering range with respect to the range of the data. This
causes the prediction to equals the mean when the data locations are
too far apart. This effect is erased in the prediction map of the
combination while preserving the accuracy in less dense areas.

The prediction standard deviations are much higher near the poles
because of the lack of data. An increased estimated uncertainty in
regions of missing data is especially visible for covariance tapering
and the combination (figs. \ref{fig:tomstaps} and \ref{fig:tomscmbs})
whereas this effect is almost imperceptible for FRK
(fig. \ref{fig:tomsfrks}). In the latter case, the predicted standard
deviation is almost constant over the globe with a very low value
compared to the two other methods. Finally, we can notice that the
predicted standard deviation obtained with the combination show some
areas with high value, mostly in the southern hemisphere. They
correspond to the locations of the selected basis functions explaining
local important departure from the mean of the TCO. They may seem a
desirable effect as they correspond to additional parameters with
respect to the covariance tapering approach.

Table \ref{tab:mspetoms} summarizes the parameter estimation results
of all three methods and shows the MSPE values obtained on the
validation dataset. Both covariance tapering and combination uses a
nugget effect with variance one while their small scale structures
differ from the weight assigned to the Wendland2 term, which is
smaller for the combination. Less than 10\% of the proposed basis
functions are selected by the LASSO. Covariance tapering and FRK show
similar results in terms of accuracy while the combination slightly
outperforms both with the lowest MSPE.

Concerning the computation times, covariance tapering is the fastest
with a running time of about 10, comprising about 9 minutes for
parameter estimation and 1 minute for prediction. Better results could
have been obtained by increasing the tapering range but the increased
memory workload caused the computing time to increase considerably
(more than one hour) due to swapping when solving the kriging
system. The combination and FRK run in approximately the same time,
around 20 minutes. Regarding the combination, the EM algorithm reaches
the stopping criterion in about two minutes with the 331 selected
basis functions within some hundredth of iterations, while 5 minutes
are necessary to build the dictionary of basis functions, without
parallelization, and 20 seconds for the lasso. Building the prediction
map takes about 2 minutes. Concerning FRK, we run the EM algorithm for
10 iterations only to reduce the computing time and obviously do not
reach the convergence criterion. This emphasizes the advantage to work
with a low number of basis functions.

\begin{table}
  \caption{Summary}
  \begin{tabular}{lccc}
    \hline\\
    model & small scale structure & number of basis functions & MSPE \\
    \hline\\

    Covariance  & nugget 4 +  & 0 & 24.63 \\
    tapering & 125Wendland2(275) & & \\
    &&&\\
    Fixed rank  & nugget 28  & 3608 &  24.09\\
          &&&\\
    Combination & nugget 4 +  & 331 & 21.64\\
     & 80Wendland2(275) & & \\
    \hline
  \end{tabular}
  \label{tab:mspetoms}
\end{table}

\section{Conclusion}

We have exposed here an original approach to interpolate large spatial
datasets that somehow reconciliates existing approaches. It benefits
from the advantages of both covariance tapering and low rank
modeling, while alleviating the computational burden associated with
the latter by using a limited number of basis functions. The model is
thus able to represent each scale of variation of the phenomenon under
study and achieves a good prediction accuracy. We have developed an
inference approach adapted to this model. The two components are
inferred almost separately through a step by step procedure. Both
inference and prediction steps are fast to compute.

In each of the three application examples, combining covariance
tapering and low rank decomposition outperformed covariance tapering
and the low rank approach alone. It is true that allowing a larger
tapering range would have improved the results of the former. However
using a too large tapering range can reduce drastically the speed of
the prediction by increasing the memory load. Therefore, our approach
can be seen as a way to improve the results given by the covariance
tapering by giving it a light low rank flavor: the examples indeed
show that adding some basis functions to the model allow to erase the
flaws obtained when using a too small taper in less sampled areas. The
examples also show the versatility of the model: the same initial set
of basis functions is used for the nested covariance model and the non
stationary Matérn covariance model while obtaining good prediction
results for each one.

Comparing computing times between different approaches for the kriging
of large datasets can be troublesome. They depend on the
implementation of each method and on the choice of the parameters
value, namely the range in covariance tapering and the number of basis
functions in fixed rank kriging. However, the application to the total
column ozone dataset show that our method provides comparable
computing times with exisiting appproaches, all implemented in
\verb?R? and run on a limited capacity computer.



\bibliographystyle{spbasic} 
\bibliography{../../biblio.bib} 

\end{document}